\documentclass[11pt]{amsart}
\pdfoutput=1

\usepackage{amscd,latexsym,amsthm,amsfonts,amssymb,amsmath,amsxtra,}
\usepackage[mathscr]{eucal}
\renewcommand\theequation{\thesection.\arabic{equation}}

\usepackage{cases}

\newcommand{\Mscr}{\mathscr{M}}

\newcommand{\ifm}{\mathrm{if}}

\newcommand{\PGL}{\mathrm{PGL}}
\newcommand{\GL}{\mathrm{GL}}

\newcommand{\Trm}{\mathrm{T}}

\newcommand{\omegabar}{\overline{\omega}}

\newcommand{\Cscr}{\mathscr{C}}

\newcommand{\Bscr}{\mathscr{B}}

\newcommand{\Krm}{\mathrm{\mathbf{K}}}

\newcommand{\Erm}{\mathrm{E}}

\newcommand{\F}{\mathrm{\mathbf{F}}}

\newcommand{\crm}{\mathrm{c}}

\newcommand{\Lscr}{\mathscr{L}}

\newcommand{\Gscr}{\mathscr{G}}

\newcommand{\Nscr}{\mathscr{N}}

\newcommand{\BA}{{\mathbb {A}}}

\newcommand{\BC}{{\mathbb {C}}}

\newcommand{\BN}{{\mathbb {N}}}

\newcommand{\BQ}{{\mathbb {Q}}}
\newcommand{\BR}{{\mathbb {R}}}

\newcommand{\BZ}{{\mathbb {Z}}}

\newcommand{\CA}{{\mathcal {A}}}
\newcommand{\CB}{{\mathcal {B}}}
\newcommand{\CC}{{\mathcal {C}}}

\newcommand{\CG}{{\mathcal {G}}}

\newcommand{\CL}{{\mathcal {L}}}

\newcommand{\CO}{{\mathcal {O}}}
\newcommand{\CP}{{\mathcal {P}}}

\newcommand{\CW}{{\mathcal {W}}}

\newcommand{\Fa}{{\mathfrak {a}}}
\newcommand{\Fb}{{\mathfrak {b}}}
\newcommand{\Fc}{{\mathfrak {c}}}

\newcommand{\Fg}{{\mathfrak {g}}}
\newcommand{\Fh}{{\mathfrak {h}}}

\newcommand{\Fl}{{\mathfrak {l}}}
\newcommand{\Fm}{{\mathfrak {m}}}
\newcommand{\Fn}{{\mathfrak {n}}}

\newcommand{\Fp}{{\mathfrak {p}}}
\newcommand{\Fq}{{\mathfrak {q}}}

\newcommand{\ScM}{{\mathscr {M}}}

\renewcommand{\Re}{{\mathrm{Re}}}

\newcommand{\bs}{\backslash}

\newtheorem{thm}{Theorem}[section]
\newtheorem{cor}[thm]{Corollary}
\newtheorem{lem}[thm]{Lemma}
\newtheorem{prop}[thm]{Proposition}

\newtheorem {ques/conj}[thm]{Question/Conjecture}

\newtheorem{defn}[thm]{Definition}
\newtheorem{rmk}[thm]{Remark}

\makeatletter

\newcommand{\Rmnum}[1]{\expandafter\@slowromancap\romannumeral #1@}
\makeatother

\begin{document}
\renewcommand{\theequation}{\arabic{equation}}
\numberwithin{equation}{section}

\title{Spectral Reciprocity and Hybrid Subconvexity Bound for triple product $L$-functions}

\author{Xinchen Miao}
\address{Mathematisches Institut\\ Endenicher Allee 60, Bonn, 53115, Germany}
\address{Max Planck Institute for Mathematics\\ Vivatsgasse 7, Bonn, 53111, Germany}
\email{miao@math.uni-bonn.de, olivermiaoxinchen@gmail.com, miao@mpim-bonn.mpg.de}

\date{December, 2024}

\subjclass[2020]{Primary 11F70, 11M41; Secondary 11F72}
\keywords{subconvexity; triple product $L$-functions; spectral reciprocity; spectral decomposition; local test vectors; local and global period integrals}
\thanks{The author was supported by ERC Advanced Grant  101054336 and Germany's Excellence Strategy grant EXC-2047/1 - 390685813.}

\begin{abstract}

Let $F$ be a number field with adele ring $\BA_F$, $\pi_1, \pi_2$ be two unitary cuspidal automorphic representations of $\PGL_2(\BA_F)$ with finite analytic conductor. We study the twisted first moment of the triple product $L$-function $L(\frac{1}{2}, \pi \otimes \pi_1 \otimes \pi_2)$ and the Hecke eigenvalues $\lambda_\pi (\Fl)$, where $\pi$ is a unitary automorphic representation of $\PGL_2(\BA_F)$ and $\Fl$ is an integral ideal coprimes with the finite analytic conductor $C(\pi \otimes \pi_1 \otimes \pi_2)$. The estimation becomes a reciprocity formula between different moments of $L$-functions. Combining with the ideas and estimations established in \cite{hmn} and \cite{subconvexity}, we study the subconvexity problem for the triple product $L$-function in the level aspect and give a new explicit hybrid subconvexity bound for $L(\frac{1}{2}, \pi \otimes \pi_1 \otimes \pi_2)$, allowing joint ramifications and conductor dropping range.

\end{abstract}

\maketitle

\tableofcontents

\section{Introduction, Background and History} \label{intro}

Subconvexity estimation is one of the most important and challenging problem in the theory of analytic number theory and $L$-functions. Let $F$ be a number field with adele ring $\BA_F$, and let $\Pi$ be an automorphic representation of a reductive group $G$. Let $L(s, \Pi)$ be the corresponding $L$-function associated to the representation $\Pi$. If $C(\Pi)$ denotes the analytic conductor of $L(s,\Pi)$, then the famous Phragmen-Lindelof principle gives the upper bound $C(\Pi)^{\frac{1}{4}+\epsilon}$ on the critical line $\Re(s)=\frac{1}{2}$. The subconvexity problem for $L(\frac{1}{2}, \Pi)$ is to establish a non-trivial upper bound of the shape as follows:
$$L(\frac{1}{2}, \Pi) \ll_{F,\epsilon} C(\Pi)^{\frac{1}{4}-\delta+\epsilon},$$
where $\delta$ is some positive absolute constant which is independent on $C(\Pi)$.

In the lower rank case $G=\GL_1, \GL_2$, the subconvexity problem has now been solved completely over a fixed general number field $F$, uniformly in all aspects ($t$, weight, spectral, level) \cite{subconvexity}. The main ingredients of the proof are integral representations, period integrals of certain $L$-functions (Ichino-Watson formula) and a spectral reciprocity relation between different families of $L$-functions, which we will mention later.

In the higher rank case for example $G=\GL_2 \times \GL_2$ or $G=\GL_2 \times \GL_2 \times \GL_2$, we are far from well-understood (See \cite{michel2} \cite{hu} \cite{michel0} \cite{michel1}  \cite{sparse}), especially in the case of hybrid subconvexity, joint ramifications and conductor dropping (See \cite{hmn}). For more history on the subconvexity problem, the interested readers may see \cite{mic} for more survey and details.

In this paper, we mainly focus on the subconvexity problem of the triple product $L$-function in the finite level aspect, which is the case $G=\GL_2 \times \GL_2 \times \GL_2$. Following \cite{hmn} \cite{subconvexity} and \cite{raphael2}, we will use the period integral approach and establish a spectral reciprocity formula for the twisted first moment of the triple product $L$-function. We sketch the rough idea as follows:

Let $\pi_1, \pi_2, \pi_3$ be three cuspidal automorphic representations of $\PGL_2(\BA_F)$ with finite analytic conductor, the Ichino-Watson triple product formula gives a rough identity as follows (See Proposition \ref{PropositionIntegralRepresentation}):
$$ \left \vert \int_{[G]} \prod_{i=1}^3 \varphi_i(g) dg\right \rvert^2= \vert \langle \varphi_3, \overline{\varphi_1 \varphi_2} \rangle \rvert^2 \asymp L(\frac{1}{2}, \pi_1 \times \pi_2 \times \pi_3) \times \prod_v I_v(\varphi_{1,v}, \varphi_{2,v}, \varphi_{3,v}).$$
Here $[G]:= Z(\BA_F)\GL_2(F) \bs \GL_2(\BA_F)$, $\varphi_i \in \pi_i$ for $i=1,2,3$ and the local integral $I_v(\varphi_{1,v}, \varphi_{2,v}, \varphi_{3,v})$ is an integral of products of local matrix coefficients for the triple product $L$-function.

Hence, we may consider the period:
$$ \int_{[G]} \varphi_1^{\Fp} \varphi_2^{\Fp} \overline{\varphi_1 \varphi_2}= \langle \varphi_1^{\Fp}\varphi_2^{\Fp}, \varphi_1 \varphi_2 \rangle,$$
where for any finite place $v$, $\varphi_i^{v}$ is the right translation of $\varphi_i$ by the matrix $\begin{pmatrix} 1 & \\ & \varpi_v \end{pmatrix}$ and $\varpi_v$ is the uniformizer for the local field $F_v$. Using the Hecke relations, the above inner product roughly equals to $\frac{1}{p^{1/2}} \sum_{\pi} \lambda_{\pi}(\Fp) \sum_{\varphi \in \CB(\pi)} \vert \langle \varphi, \varphi_1 \varphi_2 \rangle \rvert^2$, where $\CB(\pi)$ denotes an orthonormal basis in $\pi$. By the Ichino-Watson formula (See Proposition \ref{PropositionIntegralRepresentation}), the inner product in the orthonormal basis roughly equals to $L(\frac{1}{2}, \pi \otimes \pi_1 \otimes \pi_2)$. Now by the trick in the identity, we have
$$ \langle \varphi_1^{\Fp}\varphi_2^{\Fp}, \varphi_1 \varphi_2 \rangle= \langle \varphi_1^{\Fp} \overline{\varphi_1}, \overline{\varphi_2^{\Fp}} \varphi_2 \rangle.$$
By the spectral decomposition, the right hand of the above equation roughly equals to $\sum_{c(\sigma) \vert \Fp} \sum_{\psi \in \CB(\sigma, \Fp)} $ $\langle \varphi_1^{\Fp} \overline{\varphi_1}, \psi \rangle \cdot \langle \psi, \overline{\varphi_2^{\Fp}} \varphi_2 \rangle$. Now applying Ichino-Watson formula again, we obtain a spectral identity (spectral reciprocity formula) between different families of moments of $L$-functions, roughly of the following shape:
$$ \sum_{\pi} h(\pi) \cdot L(\frac{1}{2}, \pi_1 \otimes \pi_2 \otimes \pi) \approx 1+ \sum_{\sigma} \tilde{h}(\sigma) \cdot \sqrt{L(\frac{1}{2}, \pi_1 \otimes \pi_1 \otimes \sigma)L(\frac{1}{2}, \pi_2 \otimes \pi_2 \otimes \sigma)},$$
where $\pi$ and $\sigma$ run over cuspidal automorphic representations of $\PGL_2(\BA_F)$. Here we ignore the contributions from the continuous spectrum and the main term $1$ comes from the contribution from the inner product of the constant term.

In order to state our results in a more precise way, we need to give some definitions of notations.

Let $\pi_1, \pi_2$ be two unitary cuspidal automorphic representations with finite coprime conductor $\Fm$ and $\Fn$ (i.e.$(\Fm, \Fn)=1$) and bounded archimedean (spectral) parameters, $\pi_3$ be a unitary automorphic representation of $\PGL_2(\BA_F)$ with finite conductor $\Fa$. Here $\Fm,\Fn, \Fa$ are three integral ideals of $\CO_F$, where $\CO_F$ is the ring of integers of the fixed number field $F$. For all archimedean place $v \vert \infty$, we assume that $\pi_{1,v}, \pi_{2,v}, \pi_{3,v}$ are all unramified principal series representations. The norm of integral ideals $\Fm,\Fn, \Fa$ are $m,n,a$. Then we have $(m,n)=1$. Hence, the analytic conductors satisfy $C(\pi_1)=m$, $C(\pi_2)=n$, $C(\pi_3)=a$. We denote
$$ Q:=C(\pi_1 \otimes \pi_2 \otimes \pi_3)= Q_{\infty} Q_f= \prod_v Q_v= \prod_{v \vert \infty} Q_v \prod_{v< \infty} Q_v= \prod_v C_v(\pi_1 \otimes \pi_2 \otimes \pi_3),$$
and
$$ P:= P_{\infty} P_f= \prod_v P_v= \prod_{v \vert \infty} P_v  \prod_{v<\infty} P_v= \prod_v \frac{Q_v^{1/2}}{\max_{i=1,2}\{C_v(\pi_i \otimes \overline{\pi_i})\}}.$$
We note that $P_v=Q_v=1$ for almost every place $v$ and $P_v \geq 1$ at least when $Q_v$ is large enough. In the special case $(\Fa,\Fm \Fn)=1$, we have $Q=(amn)^4=a^4m^4n^4$. We let $\Fq$ be an integral ideal of $\CO_F$ with $\Fq \vert \Fa$ and $(\Fq, \Fm \Fn)=1$ ($\Fq$ maybe trivial). We denote the norm of $\Fq$ as $q$ and $q \vert a$. We assume that such integral ideal $\Fq$ is maximal, i.e. for any other $\Fq'$ satisfying the above conditions, we have $\Fq' \vert \Fq$. We write $\Fa= \Fq \Fa'$ as the ideal factorization with $(\Fq,\Fa')=1$ and we see that the integral ideal $\Fa'$ has a common factor with either $\Fm$ or $\Fn$ (The common factor may still be trivial, i.e. $\Fa'$ and $\Fm\Fn$ can be coprime).
Hence
\begin{equation*}
\begin{aligned}
  Q:&= C(\pi_1 \otimes \pi_2 \otimes \pi_3)= \prod_v C_v(\pi_1 \otimes \pi_2 \otimes \pi_3) \\
 &= \prod_{v \mid \Fq} C_v(\pi_1 \otimes \pi_2 \otimes \pi_3) \times \prod_{v \nmid \Fq} C_v(\pi_1 \otimes \pi_2 \otimes \pi_3)= q^4 \cdot \prod_{v \nmid \Fq} C_v(\pi_1 \otimes \pi_2 \otimes \pi_3).
 \end{aligned}
 \end{equation*}
We write $Q_{\Fq}:= Q_{\Fq,f} \times Q_{\infty} = \prod_{v \nmid \Fq} C_v(\pi_1 \otimes \pi_2 \otimes \pi_3)= \prod_{v \nmid \Fq< \infty} C_v(\pi_1 \otimes \pi_2 \otimes \pi_3) \cdot \prod_{v \vert \infty} C_v(\pi_1 \otimes \pi_2 \otimes \pi_3)$ and $Q=q^4 \cdot Q_{\infty} \cdot Q_{\Fq,f}$. Similarly, we can also write $P=\prod_{v \vert \infty} P_v \prod_{v \vert \Fq} P_v \prod_{v \nmid \Fq<\infty} P_v$. Since $\prod_{v \vert \Fq} P_v=q^2$, we have $P=q^2 \cdot P_{\infty} \cdot P_{\Fq,f}$.

We let the real number $\theta$ be the best exponent toward the Ramanujan-Petersson Conjecture for $\GL(2)$ over the number field $F$, we have $0 \leq \theta \leq \frac{7}{64}.$

Let $\Fl$ be an integral ideal of norm $\ell$. We assume that $(\Fl, \Fa\Fm\Fn)=1$, hence $(l, amn)=1$. We define the following for the cuspidal contribution:

\begin{equation}\label{CuspidalPart}
\mathscr{C}(\pi_1,\pi_2,\Fa,\Fm, \Fn, \Fl):= C_1\cdot \sum_{\substack{\pi \ \mathrm{cuspidal} \\ C(\pi) \vert \Fc \cdot [\Fm,\Fn, \Fa]}}\lambda_\pi(\Fl)\frac{L(\frac{1}{2}, \pi\otimes\pi_1\otimes\pi_2)}{\Lambda(1, \pi,\mathrm{Ad})}f(\pi_\infty)H(\pi,\Fa, \Fm, \Fn),
\end{equation},

where $[\Fm,\Fn,\Fa]$ means the least common multiple of the three integral ideals ($[\Fm,\Fn,\Fa]$ is the minimal integral ideal satisfying $\Fm,\Fn,\Fa \vert [\Fm,\Fn,\Fa]$, and for any integral ideal $\Fg$ satisfying $\Fm,\Fn,\Fa \vert \Fg$, we have $[\Fm,\Fn,\Fa] \vert \Fg$) and $\Fc$ is a fixed integral ideal with its norm a fixed positive absolute bounded integer. The integral ideal $\Fc$ is coprime with the integral ideal $\Fq$ and is determined by the choices of test vectors (See Section \ref{SectionEstimation} and also Section 6.3 Choice of test vectors and Proposition 6.5 in \cite{hmn}). Moreover, $H$ is certain weight function in terms of finite many ramified non-archimedean local places defined in Section \ref{Connection} and $f(\pi_\infty)$ is defined in Section \ref{SectionInterlude}. Here the constant $C_1$ is a positive constant depending only on the number field $F$ and the nature of the three representations $\pi, \pi_1, \pi_2$. If $\pi,\pi_1,\pi_2$ are all cuspidal, then $C_1=2 \Lambda_F(2)$. If we further assume that the characteristic of all the residue fields corresponding to the prime ideals in $\Fa,\Fm,\Fn$ are large enough, then we can take $\Fc=\textbf{1}$ trivially. If the characteristic of the residue fields are bounded, then the integral ideal $\Fc$ may not be trivial, however, it is fixed and the corresponding norm is a fixed non-negative absolute bounded integer
(See Section 6.3 Choice of test vectors and Proposition 6.5 in \cite{hmn}). 

For the continuous part, we denote by $\pi_{\omega}(it)$ the principal series $\omega|\cdot|^{it}\boxplus \omegabar|\cdot|^{-it}$ and define similarly
\begin{equation}\label{ContinuousPart}
\begin{split}
\mathscr{E}(\pi_1,\pi_2,\Fa,\Fm,\Fn, \Fl):= C_2\cdot \sum_{\substack{\omega\in\widehat{\F^\times\setminus \BA_\F^{1}} \\ C(\omega)^2 \vert \Fc \cdot [\Fm,  \Fn,\Fa]}}&\int_{-\infty}^\infty \lambda_{\pi_\omega(it)}(\Fl)f(\pi_{\omega_\infty}(it))H(\pi_\omega(it),\Fa,\Fm,\Fn) \\ \times & \frac{L(\tfrac{1}{2}+it, \pi_1\otimes\pi_2\otimes\omega)L(\tfrac{1}{2}-it, \pi_1\otimes\pi_2\otimes\omegabar)}{\Lambda^*(1, \pi_\omega(it),\mathrm{Ad})} \frac{dt}{4\pi}.
\end{split}
\end{equation}
In this case, $C_2= 2 \Lambda_F^{*}(1)$ and $H$ is certain weight function in terms of finite many ramified non-archimedean local places defined in Section \ref{Connection}. We also note that the completed $L$-functions satisfy $\Lambda(s,\pi,\mathrm{Ad})=\Lambda(s,\chi^2)\Lambda(s,\chi^{-2})\zeta_F(s)$, where $\pi$ is an Eisenstein series normalized induced from a character $\chi$. In above Equation \ref{ContinuousPart}, $\chi=\omega \vert \cdot \rvert^{it}$ and $\omega$ is a unitary Hecke character. For $\chi^2 \neq 1$, we define $\Lambda^{*}(1,\pi,\mathrm{Ad})=\Lambda(1,\chi^2)\Lambda(1,\chi^{-2})\zeta_F^{*}(1)$, where $\zeta_F^{*}(1)$ is the residue of the Dedekind zeta function at $s=1$, and is a positive real number by the class number formula. We also note that the Dedekind zeta function has a simple pole at $s=1$. When $\chi^2=1$, we define $1/ \Lambda^{*}(1,\pi,\mathrm{Ad}):=0$ (Section 3.2 in \cite{BJN}). Hence, the function $1/ \Lambda^{*}(1,\pi,\mathrm{Ad})$ is continuous in terms of the induced character $\chi$. Since the Dedekind zeta function has a simple pole at $s=1$, if $\chi^2=1$ and $\pi$ is normalized induced from $\chi$, for some small real number $t$ satisfying $\vert t \rvert \leq 1$ (can take zero), we have $1/ \Lambda^{*}(1+it,\pi,\mathrm{Ad}) \gg_F t^2$.

We define
\begin{equation}\label{DefinitionMoment1}
\mathscr{M}(\pi_1,\pi_2,\Fa,\Fm,\Fn,\Fl) :=\mathscr{C}(\pi_1,\pi_2,\Fa,\Fm,\Fn,\Fl)+ \mathscr{E}(\pi_1,\pi_2,\Fa,\Fm,\Fn,\Fl).
\end{equation}
The first theorem establishes an upper bound for this twisted first moment.

\begin{thm} \label{moment}
Let $\pi_1,\pi_2$ be two unitary $\theta_i$-tempered ($i=1,2$) cuspidal automorphic representations with bounded archimedean (spectral) parameters and finite coprime conductor $\mathfrak {m}$ and $\mathfrak {n}$. We let the real number $\theta_i$ be the best exponent toward the Ramanujan-Petersson Conjecture for $\GL(2)$ over the number field $F$ for $\pi_1$ and $\pi_2$, we have $0 \leq \theta_i \leq \frac{7}{64}.$
Assume that for all archimedean places $v \vert \infty$, both $\pi_{1,v}$ and $\pi_{2,v}$ are unramified principal series representation. Let $\Fq, \Fl$ be two coprime ideals of $\CO_F$ with the condition $(\Fq \Fl, \mathfrak {m} \mathfrak {n})=1$ and write $q$ and $\ell$ for their respective norms. Then the twisted first moment satisfies
\begin{equation}  \label{moment1}
 \ScM(\pi_1, \pi_2, \Fa,\Fm,\Fn, \Fl) \ll_{\pi_{1,\infty}, \pi_{2,\infty}, F, \varepsilon} (mna\ell)^{\epsilon} \cdot \left(\ell^{\frac{3}{2}+2\theta_1+2\theta_2} \cdot P_{\Fq,f}^{-1/4+\frac{\theta}{2}} \cdot q^{-1/2+\theta}+  \ell^{-1/2+\theta_1+ \theta_2} \right).
\end{equation}

\end{thm}

\begin{cor}
With all the notations same as above, we have the following twisted first moment estimation:
\begin{equation}  \label{moment2}
 \ScM(\pi_1, \pi_2, \Fa,\Fm,\Fn, \Fl) \ll_{\pi_{1,\infty}, \pi_{2,\infty}, F, \varepsilon} (mna \ell)^{\epsilon} \cdot \left(\ell^{\frac{3}{2}+2\theta_1+2\theta_2} \cdot P_{f}^{-1/4+\frac{\theta}{2}} +  \ell^{-1/2+\theta_1+\theta_2} \right).
\end{equation}
\end{cor}

Combining Theorem \ref{moment} with the amplification method, we obtain the following subconvexity bounds in the hybrid level aspect.

\begin{thm} \label{subconvex1}
Let $\pi_1,\pi_2$ be two unitary cuspidal automorphic representations and $\pi_3$ be unitary automorphic representation with corresponding finite levels defined previously. Assume that for all archimedean places $v \vert \infty$, $\pi_{1,v}$, $\pi_{2,v}$ and $\pi_{3,v}$ are all unramified principal series representation. We have the following subconvex estimation:
\begin{equation}\label{SubConv1}
L\left( \tfrac{1}{2}, \pi_1 \otimes\pi_2 \otimes\pi_3 \right) \ll_{\varepsilon, F, \pi_{1,\infty},\pi_{2,\infty},\pi_{3,\infty}} Q_f^{1/4+\varepsilon} \cdot P_f^{-(\frac{1}{4}-\frac{\theta}{2})(1-2\theta_1-2\theta_2)/(7-2\theta_1-2\theta_2)}.
\end{equation}
If we pick $\theta=\theta_1=\theta_2=\frac{7}{64}$, then we have $(\frac{1}{4}-\frac{\theta}{2})(1-2\theta_1-2\theta_2)/(7-2\theta_1-2\theta_2)>1/60$. Hence, we have
$$L\left( \tfrac{1}{2}, \pi_1 \otimes\pi_2 \otimes\pi_3\right) \ll_{\varepsilon, F, \pi_{1,\infty},\pi_{2,\infty},\pi_{3,\infty}} Q_f^{1/4+\varepsilon} \cdot P_f^{-1/60}$$
unconditionally.
\end{thm}

\begin{rmk}
This is an explicit version for the hybrid subconvexity in the level aspect for the triple product $L$-function (See Theorem 1.3 in \cite{hmn}). Assume that all the archimedean components of three representations $\pi_1,\pi_2,\pi_3$ are unramified principal series, we can pick the absolute constant $\delta=\frac{15}{896}>\frac{1}{60}$. If we further assume the Ramanujan-Petersson conjecture, we can have $\delta=\frac{1}{28}$.
\end{rmk}

In the special case that all the three levels are coprime to each other, we have a more explicit hybrid subconvexity bound for the triple product $L$-functions in the level (also depth) aspect which states as follows:

\begin{thm} \label{subconvex1'}
 Let $\pi_1,\pi_2, \pi_3$ be three unitary cuspidal automorphic representations with finite conductor $\mathfrak {m}, \mathfrak{n}$ and $\mathfrak {a}$ which are coprime to each other. Assume that for all archimedean places $v \vert \infty$, $\pi_{1,v}$, $\pi_{2,v}$ and $\pi_{3,v}$ are all unramified principal series representation. We further assume that the integral ideals $\mathfrak m$ and $\mathfrak n$ are squarefull (squarefull ideals are integral ideals for which all the prime ideal factors exponents are at least two). We write the prime ideal factorizations $\mathfrak m= \prod_{i=1}^s \mathfrak m_i^{u_i}$ and $\mathfrak n= \prod_{i=1}^t \mathfrak n_i^{v_i}$, where $\mathfrak m_i, \mathfrak n_i$ are all coprime prime ideals with norm $m_i, n_i$ and $u_i, v_i \geq 2$. We denote $M:= \prod_{i=1}^s m_i$ and $N:= \prod_{i=1}^t n_i$. If $(\Fa, \mathfrak {m} \mathfrak {n})=1$, then for any $\epsilon>0$, we have the following subconvex estimation:
\begin{equation}\label{SubConv1'}
\begin{aligned}
 L\left( \tfrac{1}{2}, \pi_1 \otimes\pi_2 \otimes\pi_3 \right) & \ll_{\varepsilon, F, \pi_{1,\infty},\pi_{2,\infty},\pi_{3,\infty}}  (MN)^{\frac{(1/4-\theta/2)(1-2\theta_1-2\theta_2)}{7-2\theta_1-2\theta_2}} \\ & \times m^{1-\frac{(1/4-\theta/2)(1-2\theta_1-2\theta_2)}{7-2\theta_1-2\theta_2}+\epsilon}  n^{1-\frac{(1/4-\theta/2)(1-2\theta_1-2\theta_2)}{7-2\theta_1-2\theta_2}+\epsilon}  a^{1-\frac{(1/2-\theta)(1-2\theta_1-2\theta_2)}{7-2\theta_1-2\theta_2}+\epsilon}.
\end{aligned}
\end{equation}
\end{thm}

By picking $a=\max\{m,n,a\}$, we have the following corollary.

\begin{cor}
By using the same notations in Theorem \ref{subconvex1'}, we have
\begin{equation}
 L\left( \tfrac{1}{2}, \pi_1 \otimes\pi_2 \otimes\pi_3 \right)  \ll_{\varepsilon, F, \pi_{1,\infty},\pi_{2,\infty},\pi_{3,\infty}}  (MN)^{\frac{(1/4-\theta/2)(1-2\theta_1-2\theta_2)}{7-2\theta_1-2\theta_2}}  \times (mna)^{1-\frac{(1-2\theta)(1-2\theta_1-2\theta_2)}{21-6\theta_1-6\theta_2}+\epsilon},
\end{equation}
and
\begin{equation}
 L\left( \tfrac{1}{2}, \pi_1 \otimes\pi_2 \otimes\pi_3 \right)  \ll_{\varepsilon, F, \pi_{1,\infty},\pi_{2,\infty},\pi_{3,\infty}}  (mna)^{1-\frac{(5/2-5\theta)(1-2\theta_1-2\theta_2)}{63-18\theta_1-18\theta_2}}.
\end{equation}
If we do not assume that two integral ideals $\mathfrak{m}$ and $\mathfrak{n}$ are squarefull, we will have
\begin{equation}
 L\left( \tfrac{1}{2}, \pi_1 \otimes\pi_2 \otimes\pi_3 \right)  \ll_{\varepsilon, F, \pi_{1,\infty},\pi_{2,\infty},\pi_{3,\infty}}  (mna)^{1-\frac{(1/2-\theta)(1-2\theta_1-2\theta_2)}{21-6\theta_1-6\theta_2}}.
\end{equation}
\end{cor}

If the automorphic representation $\pi_3$ is an Eisenstein series, we have

\begin{thm} \label{subconvex2}
Let $\pi_1,\pi_2$ be two unitary cuspidal automorphic representations with finite coprime conductor $\mathfrak {m}$ and $\mathfrak {n}$ which are defined in Theorem \ref{subconvex1'}. Let $\Fh$ be an integral ideal of norm $h$ and $\chi$ a unitary Hecke character with finite conductor $\Fh$. All the integral ideals listed here are coprime to each other. We further assume that the integral ideals $\mathfrak m$ and $\mathfrak n$ are squarefull (squarefull ideals are integral ideals for which all the prime ideal factors exponents are at least two). We write the prime ideal factorizations $\mathfrak m= \prod_{i=1}^s \mathfrak m_i^{u_i}$ and $\mathfrak n= \prod_{i=1}^t \mathfrak n_i^{v_i}$, where $\mathfrak m_i, \mathfrak n_i$ are all coprime prime ideals with norm $m_i, n_i$ and $u_i, v_i \geq 2$. We denote $M:= \prod_{i=1}^s m_i$ and $N:= \prod_{i=1}^t n_i$. If $(\Fh, \mathfrak {m} \mathfrak {n})=1$, then for any $\epsilon>0$, we have the following subconvex estimation:
\begin{equation}\label{SubConv2}
\begin{aligned}
 L\left( \tfrac{1}{2}, \pi_1 \otimes\pi_2 \otimes \chi \right) & \ll_{\varepsilon, F, \pi_{1,\infty},\pi_{2,\infty},\chi_\infty}  (MN)^{\frac{(1/4-\theta/2)(1-2\theta_1-2\theta_2)}{7-2\theta_1-2\theta_2}} \\ & \times m^{1-\frac{(1/4-\theta/2)(1-2\theta_1-2\theta_2)}{7-2\theta_1-2\theta_2}+\epsilon}  n^{1-\frac{(1/4-\theta/2)(1-2\theta_1-2\theta_2)}{7-2\theta_1-2\theta_2}+\epsilon}  h^{1-\frac{(1/2-\theta)(1-2\theta_1-2\theta_2)}{7-2\theta_1-2\theta_2}+\epsilon}.
\end{aligned}
\end{equation}
Similarly, we will also have
\begin{equation}
 L\left( \tfrac{1}{2}, \pi_1 \otimes\pi_2 \otimes\chi \right)  \ll_{\varepsilon, F, \pi_{1,\infty},\pi_{2,\infty},\pi_{3,\infty}}  (MN)^{\frac{(1/4-\theta/2)(1-2\theta_1-2\theta_2)}{7-2\theta_1-2\theta_2}}  \times (mnh)^{1-\frac{(1-2\theta)(1-2\theta_1-2\theta_2)}{21-6\theta_1-6\theta_2}+\epsilon},
\end{equation}
and
\begin{equation}
 L\left( \tfrac{1}{2}, \pi_1 \otimes\pi_2 \otimes\chi \right)  \ll_{\varepsilon, F, \pi_{1,\infty},\pi_{2,\infty},\pi_{3,\infty}}  (mnh)^{1-\frac{(5/2-5\theta)(1-2\theta_1-2\theta_2)}{63-18\theta_1-18\theta_2}}.
\end{equation}
If we do not assume that two integral ideals $\mathfrak{m}$ and $\mathfrak{n}$ are squarefull, we will have
\begin{equation}
 L\left( \tfrac{1}{2}, \pi_1 \otimes\pi_2 \otimes\chi \right)  \ll_{\varepsilon, F, \pi_{1,\infty},\pi_{2,\infty},\pi_{3,\infty}}  (mnh)^{1-\frac{(1/2-\theta)(1-2\theta_1-2\theta_2)}{21-6\theta_1-6\theta_2}}.
\end{equation}
\end{thm}

The results here are even new when the ground field $\F=\BQ$.

\section{Automorphic Forms Preliminaries} \label{pre}

In this paper, $F/\mathbb{Q}$ will denote a fixed number field with ring of intergers $\CO_F$ and discriminant $\Delta_F$. We make the assumption that all prime ideals considering in this paper ($\Fq,\Fl, \Fa,\mathfrak{m},\mathfrak{n},\Fc, \Fh$) do not divide $\Delta_F$. We let $\Lambda_F$ be the complete Dedekind $\zeta$-function of $F$; it has a simple pole at $s=1$ with residue $\Lambda_F^*(1)$.

For $v$ a place of $F$, we set $F_v$ for the completion of $F$ at the place $v$. We will also write $F_{\Fp}$ if $v$ is finite place that corresponds to a prime ideal $\Fp$ of $\CO_F$. If $v$ is non-Archimedean, we write $\CO_{F_v}$ for the ring of integers in $F_v$ with maximal ideal $\Fm_v$ and uniformizer $\varpi_v$. The size of the residue field is $q_v=\CO_{F_v}/\Fm_v$. For $s\in\BC$, we define the local zeta function $\zeta_{F_v}(s)$ to be $(1-q_v^{-s})^{-1}$ if $v<\infty$, $\zeta_{F_v}(s)=\pi^{-s/2}\Gamma(s/2)$ if $v$ is real and $\zeta_{F_v}(s)=2(2\pi)^{-s}\Gamma(s)$ if $v$ is complex.

The adele ring of $F$ is denoted by $\BA_F$ and its unit group $\BA^\times_F$. We also set $\widehat{\CO}_F:=\prod_{v<\infty} \CO_{F_v}$ for the profinite completion of $\CO_F$ and $\BA^1_F=\{ x\in \BA_F^\times \; : \; \vert x \rvert=1 \}$, where $\vert \cdot \rvert : \BA_F^\times \rightarrow \BR_{>0}$ is the adelic norm map.

We denote by $\psi = \prod_v \psi_v$ the additive character $\psi_{\BQ} \circ \text{Tr}_{F/ \BQ}$ where $\psi_{\BQ}$ is the additive character on $\BQ \setminus \BA_{\BQ}$ with value $e^{2\pi i x}$ on $\BR$. For $v<\infty$, we let $d_v$ be the conductor of $\psi_v$, which is the smallest non-negative integer such that $\psi_v$ is trivial on $\Fm_v^{d_v}$. In this case, we have $\Delta_F=\prod_{v<\infty} q_v^{d_v}$. We also set $d_v=0$ for the Archimedean local place $v$.

If $R$ is a commutative ring, $\GL_2(R)$ is by definition the group of $2\times 2$ matrices with coefficients in $R$ and determinant in the multiplicative group $R^{\times}$. We also define the following standard subgroups:
$$B(R)=\left\{\begin{pmatrix} a & b \\ & d\end{pmatrix} \; : \; a,d \in R^{\times}, b\in R\right\}, \; P(R)= \left\{ \begin{pmatrix} a & b \\ & 1\end{pmatrix} \; : \; a\in R^{\times}, b\in R\right\}, $$
$$Z(R)=\left\{\begin{pmatrix} z & \\ & z\end{pmatrix} \; : \; z\in R^{\times}\right\}, \; A(R)=\left\{\begin{pmatrix} a &  \\ & 1\end{pmatrix} \; : \; a\in R^{\times}\right\}, $$ $$N(R)=\left\{\begin{pmatrix} 1 & b \\ & 1\end{pmatrix} \; : \; b\in R\right\}.$$
We also set
$$n(x)=\begin{pmatrix}
1 & x \\ & 1
\end{pmatrix}, \hspace{0.4cm}w = \begin{pmatrix} & 1 \\ -1 & \end{pmatrix} \hspace{0.4cm} \mathrm{and} \hspace{0.4cm} a(y)= \begin{pmatrix}
y & \\ & 1
\end{pmatrix}.
$$
For any place $v$, we let $K_v$ be the maximal compact subgroup of $G(F_v)$ defined by
\begin{equation}\label{Compact}
K_v= \left\{ \begin{array}{lcl}
\GL_2(\CO_{F_v}) & \text{if} & v \; \mathrm{is \; finite} \\
 & & \\
\mathrm{O}_2(\BR) & \text{if} & v \; \mathrm{is \; real} \\
 & & \\
\mathrm{U}_2(\BC) & \text{if} & v \; \mathrm{is \; complex}.
\end{array}\right.
\end{equation}

We also set $K:= \prod_v K_v$. If $v<\infty$ and $n\geqslant 0$, we define 
$$K_{v,0}(\varpi_v^n):= \left\{ \begin{pmatrix}
a & b \\ c & d
\end{pmatrix} \in K_v \; : \;  c \in \Fm_v^n\right\}.$$
If $\Fb$ is an integral ideal of $\CO_F$ with prime factorization $\Fb=\prod_{v<\infty}\Fp_v^{f_v(\Fb)}$ ($\Fp_v$ is the prime ideal corresponding to the finite place $v$), then we set
$$K_0(\Fb):=\prod_{v<\infty} K_{v,0}\left(\varpi_v^{f_v(\Fb)} \right).$$

We use the same measures normalizations as in \cite{subconvexity}. At each place $v$, $dx_v$ denotes a self-dual measure on $F_v$ with respect to $\psi_v$. If $v<\infty$, $dx_v$ gives the measure $q_v^{-d_v/2}$ to $\CO_{F_v}$. We define $dx=\prod_v dx_v$ on $\BA_F$. We take $d^\times x_v=\zeta_{F_v}(1)\frac{dx_v}{\vert x_v \rvert}$ as the Haar measure on the multiplicative group $F_v^\times$ and $d^\times x = \prod_v d^\times x_v$ as the Haar measure on the idele group $\BA^\times_F$.
We provide $K_v$ with the probability Haar measure $dk_v$. We identify the subgroups $Z(F_v)$, $N(F_v)$ and $A(F_v)$ with respectively $F_v^\times,$ $F_v$ and $F_v^\times$ and equipped them with the measure $d^\times z$, $dx_v$ and $d^\times y_v$. Using the Iwasawa decomposition, namely $\GL_2(F_v)=Z(F_v)N(F_v)A(F_v)K_v$, a Haar measure on $\GL_2(\F_v)$ is given by 
\begin{equation}\label{HaarMeasure}
dg_v = d^\times z dx_v \frac{d^\times y_v}{\vert y_v \rvert}dk_v.
\end{equation} 
The measure on the adelic points of the various subgroups are just the product of the local measures defined above. We also denote by $dg$ the quotient measure on $$X:= Z(\BA_F)\GL_2(F) \setminus \GL_2(\BA_F),$$ 
with total mass $V_F:=\mathrm{vol}(X)<\infty$. 

Let $\pi=\otimes_v\pi_v$ be a unitary automorphic representation of $\PGL_2(\BA_F)$ and fix $\psi$ a character of $\F \setminus \BA_F$. The intertwiner
\begin{equation}\label{NaturalIntertwiner}
\pi \ni \varphi \longmapsto  W_\varphi(g):=\int_{F\setminus \BA_F} \varphi(n(x)g)\psi(-x)dx,
\end{equation}
gives a $\GL_2(\BA_F)$-equivariant embedding of $\pi$ into a space of functions $W : \GL_2(\BA_F) \rightarrow \BC$ satisfying $W(n(x)g))=\psi(x)W(g)$. The image is called the Whittaker model of $\pi$ with respect to $\psi$ and it is denoted by $\CW(\pi,\psi)$. This space has a factorization $\otimes_v \CW(\pi_v,\psi_v)$ into local Whittaker models of $\pi_v$. A pure tensor $\otimes_v \varphi_v$ has a corresponding decomposition $\prod_v W_{\varphi_v}$ where $W_{\varphi_v}(1)=1$ and is $K_v$-invariant for almost all place $v$.

We define a normalized inner product on the space $\CW(\pi_v,\psi_v)$ by the rule 
\begin{equation}\label{NormalizedInnerProduct}
\vartheta_v(W_v,W_v') :=\zeta_{F_v}(2) \times \frac{\int_{F_v^\times}W_v(a(y))\overline{W}_v'(a(y))d^\times y}{\zeta_{F_v}(1)L(1, \pi_v,\mathrm{Ad})}.
\end{equation}
This normalization has the good property that $\vartheta_v(W_v,W_v)=1$ for $\pi_v$ and $\psi_v$ unramified and $W_v(1)=1$ \cite[Proposition 2.3]{classification}. We also fix for each place $v$ an invariant inner product $\langle \cdot,\cdot\rangle_v$ on $\pi_v$ and an equivariant isometry $\pi_v \rightarrow \CW(\pi_v,\psi_v)$ with respect to \eqref{NormalizedInnerProduct}.

Let $L^2(X)$ be the Hilbert space of square integrable functions $\varphi : X \rightarrow \BC$.  If $\pi$ is a cuspidal representation, for any $\varphi\in\pi$, we can define the $L^2$-norm by
\begin{equation} \label{L^2normCuspidal}
||\varphi||_{L^2}^2:= \int_{X} |\varphi(g)|^2 dg.
\end{equation}
We denote by $L_{\mathrm{cusp}}^2(X)$ the closed subspace of cusp forms, i.e. the functions $\varphi\in L^2(X)$ with the additional property that 
$$\int_{F \setminus \BA_F}\varphi(n(x)g)dg=0, \ \ \mathrm{a.e.} \ g\in \GL_2(\BA_F).$$
Each $\varphi\in L^2_{\mathrm{cusp}}(X)$ admits a Fourier expansion
\begin{equation}\label{FourierSeries}
\varphi(g)= \sum_{\alpha \in F^\times} W_\varphi \left(\begin{pmatrix} \alpha & \\ & 1 \end{pmatrix} g \right),
\end{equation}
\begin{equation}\label{Whittaker-Cuspidal}
W_\varphi(g)=\int_{ F \setminus \BA_F}\varphi\left( \begin{pmatrix} 1 & x \\ & 1 \end{pmatrix} g \right) \psi(-x)dx.
\end{equation}
The group $\GL_2(\BA_F)$ acts by right translations on both spaces $L^2(X)$ and $L_{\mathrm{cusp}}^2(X)$ and the resulting representation is unitary with respect to \eqref{L^2normCuspidal}. It is well known that each irreducible component $\pi$ decomposes into $\pi = \otimes_v \pi_v$ where $\pi_v$ are smooth irreducible and unitary representations of the local group $\GL_2(F_v)$. The spectral decomposition is established in the first four chapters of \cite{analytic} and gives the orthogonal decomposition
\begin{equation}\label{OrthogonalDecomposition}
L^2(X)=L^2_{\mathrm{cusp}}(X)\oplus L^2_{\mathrm{res}}(X)\oplus L^2_{\mathrm{cont}}(X).
\end{equation}
$L^2_{\mathrm{cusp}}(X)$ decomposes as a direct sum of irreducible $\GL_2(\BA_F)$-representations which are called the cuspidal automorphic representations. $L^2_{\mathrm{res}}(X)$ is the sum of all one dimensional subrepresentations of $L^2(X)$. Finally the continuous part $L^2_{\mathrm{cont}}(X)$ is a direct integral of irreducible $\GL_2(\BA_F)$-representations and it is expressed via the Eisenstein series. In this paper, we call the irreducible components of $L^2_{\mathrm{cusp}}$ and $L^2_{\mathrm{cont}}$ the unitary automorphic representations. If $\pi$ is a unitary representation appearing in the continuous part, we say that $\pi$ is Eisenstein.

For any ideal $\Fb$ of $\CO_F$, we write $L^2(X,\Fb):= L^2(X)^{K_0(\Fb)}$ for the subspace of level $\Fb$ automorphic forms, which is the closed subspace of functions that are invariant under the subgroup $K_0(\Fb)$.

Note that if $\pi$ is a cuspidal representation, we have a unitary structure on $\pi$ given by \eqref{L^2normCuspidal}. If $\pi$ belongs to the continuous spectrum and $\varphi$ is the Eisenstein series associated to a section $f : \GL_2(\BA_F) \rightarrow \BC$ in an induced representation of $B(\BA_F)$ (see for example \cite[Section 4.1.6]{subconvexity} for the basic facts and notations concerning Eisenstein series), we can define the norm of $\varphi$ by setting 
$$||\varphi||^2_{\mathrm{Eis}}:= \int_{K}|f(k)|^2dk.$$
We define the canonical norm of $\varphi$ by
\begin{equation}\label{CanonicalNorm}
||\varphi||^2_{\mathrm{can}} := \left\{ \begin{array}{lcl}
||\varphi||_{L^2(X)}^2 & \text{if} & \pi \; \mathrm{is \; cuspidal} \\
 & & \\
2\Lambda_\F^*(1) ||\varphi||_{\mathrm{Eis}}^2 & \text{if} & \pi \; \mathrm{is \; Eisenstein},
\end{array}\right.
\end{equation}
Using \cite[Lemma 2.2.3]{subconvexity}, we can compare the global and the local inner product : for $\varphi=\otimes_v \varphi_v \in \pi=\otimes_v\pi_v$ a pure tensor with $\pi$ either cuspidal or Eisenstein and non-singular, i.e. $\pi=\chi_1\boxplus\chi_2$ with $\chi_i$ unitary, $\chi_1\chi_2=1$ and $\chi_1\neq\chi_2,$ we have
\begin{equation}\label{Comparition}
||\varphi||_{\mathrm{can}}^2=2 \Delta_\F^{1/2} \Lambda^*(1,\pi,\mathrm{Ad})\prod_v \langle \varphi_{v},\varphi_v\rangle_v,
\end{equation}
where $\Lambda(s,\pi,\mathrm{Ad})$ is the complete adjoint $L$-function $\prod_v L(s,\pi,\mathrm{Ad})$ and $\Lambda^*(1,\pi,\mathrm{Ad})$ is the first nonvanishing coefficient in the Laurent expansion around $s=1$.
This regularized value satisfies \cite{adjoint}
\begin{equation}\label{BoundAdjoint}
\Lambda^*(1,\pi,\mathrm{Ad})=\mathrm{C}(\pi)^{o(1)}, \; \; \mathrm{as} \; \mathrm{C}(\pi)\rightarrow \infty,
\end{equation}
where $\mathrm{C}(\pi)$ is the analytic conductor of $\pi$, as defined in \cite[Section 1.1]{subconvexity}.

\section{Integral representations of triple product $L$-functions}\label{SectionRankin} 

Let $\pi_1,\pi_2,\pi_3$ be three unitary automorphic representations of $\PGL_2(\BA_F)$ such that at least one of them is cuspidal, say $\pi_2$. We consider the linear functional on $\pi_1\otimes\pi_2\otimes\pi_3$ defined by
$$I (\varphi_1\otimes\varphi_2\otimes \varphi_3):= \int_X \varphi_1(g)\varphi_2(g)\varphi_3(g)dg.$$
This period is closely related to the central value of the triple product $L$-function $L(\tfrac{1}{2}, \pi_1\otimes\pi_2\otimes\pi_3)$. In order to state the result, we write $\pi_i=\otimes_v \pi_{i,v}$ and for each $v$, we can consider the matrix coefficient
\begin{equation}\label{DefinitionMatrixCoefficient}
I'_v(\varphi_{1,v}\otimes\varphi_{2,v}\otimes\varphi_{3,v}) :=\int_{\PGL_2(\F_v)}\prod_{i=1}^3\langle \pi_{i,v}(g_v)\varphi_{i,v},\varphi_{i,v}\rangle_v dg_v. 
\end{equation}
It is a fact that \cite[(3.27)]{subconvexity} 
\begin{equation}\label{Fact}
\frac{I'(\varphi_{1,v}\otimes\varphi_{2,v}\otimes\varphi_{3,v})}{\prod_{i=1}^3 \langle \varphi_{i,v},\varphi_{i,v}\rangle_v}= \zeta_{F_v}(2)^2 \frac{L(\tfrac{1}{2}, \pi_{1,v}\otimes\pi_{2,v}\otimes\pi_{3,v})}{\prod_{i=1}^3 L(1,\pi_{i,v},\mathrm{Ad})},
\end{equation}
when $v$ is non-Archimedean and all vectors are unramified. It is therefore natural to consider the normalized version
\begin{equation}\label{DefinitionNormalizedMatrixCoefficient}
I_v(\varphi_{1,v}\otimes \varphi_{2,v}\otimes \varphi_{3,v}) := \zeta_{\F_v}(2)^{-2} \frac{\prod_{i=1}^3 L(1,\pi_{i,v},\mathrm{Ad})}{L(\tfrac{1}{2}, \pi_{1,v}\otimes\pi_{2,v}\otimes\pi_{3,v})} I'_v (\varphi_{1,v}\otimes \varphi_{2,v} \otimes \varphi_{3,v}).
\end{equation}
The following proposition connects the global trilinear form $I$ with the central value $L(\tfrac{1}{2}, \pi_1\otimes\pi_2\otimes\pi_3)$ and the local matrix coefficients $I_v$. The proof when at least one of the $\pi_i$'s is Eisenstein can be found in \cite[Equation 4.21]{subconvexity} and is a consequence of the Rankin-Selberg method. The result when all $\pi_i$ are cuspidal is due to Ichino \cite{ichino}.

\begin{prop}\label{PropositionIntegralRepresentation} 
Let $\pi_1,\pi_2,\pi_3$ be unitary automorphic representations of $\PGL_2(\BA_F)$ such that at least one of them is cuspidal. Let $\varphi_i = \otimes_v \varphi_{i,v}\in \otimes_v \pi_{i,v}$ be pure tensors and set $\varphi :=\varphi_1\otimes\varphi_2\otimes\varphi_3$.

\begin{enumerate}
\item If none of the $\pi_i$'s ($i=1,2,3$) is a singular Eisenstein series, then
$$
\frac{|I(\varphi)|^2}{\prod_{i=1}^3 ||\varphi_i||^2_{\mathrm{can}}} = \frac{C}{8\Delta_F^{3/2}}\cdot\frac{\Lambda(\tfrac{1}{2}, \pi_1\otimes\pi_2\otimes\pi_3)}{\prod_{i=1}^3 \Lambda^*(1, \pi_i,\mathrm{Ad})}\prod_v \frac{I_v(\varphi_v)}{\prod_{i=1}^3\langle \varphi_{i,v},\varphi_{i,v}\rangle_v},
$$
with $C=\Lambda_F(2)$ if all $\pi_i$ are cuspidal and $C=1$ if at least one $\pi_i$ is Eisenstein and non-singular.

\item Assume that $\pi_3=1 \boxplus 1$ and let $\varphi_3$ be the Eisenstein associated to the section $f_3(0) \in 1 \boxplus 1$ which for $\Re(s)>0$, is defined as follows:
$$f_3(g,s):= \vert \det(g) \rvert^s \cdot \int_{\BA_F^{\times}} \Phi((0,t)g) \vert t \rvert^{1+2s} d^{\times} t \in \vert \cdot \rvert^s \boxplus \vert \cdot \rvert^{-s},$$
where $\Phi=\prod_v \Phi_v$ and $\Phi_v=1_{\CO_{F_v}}^2$ for finite $v$. Then we have
$$
\frac{|I(\varphi)|^2}{\prod_{i=1}^2 ||\varphi_i||^2_{\mathrm{can}}} = \frac{1}{4 \Delta_F}\cdot\frac{\Lambda(\tfrac{1}{2}, \pi_1\otimes\pi_2\otimes\pi_3)}{\prod_{i=1}^2 \Lambda^*(1, \pi_i,\mathrm{Ad})}\prod_v \frac{I_v(\varphi_v)}{\prod_{i=1}^3\langle \varphi_{i,v},\varphi_{i,v}\rangle_v}.
$$
\end{enumerate}

\end{prop}

\subsection{Hecke operators} \noindent Let $\Fp$ be a prime ideal of $\CO_F$ of norm $p$ and $n\in\mathbb{N}$. Let $F_p$ be the completion of the number field $F$ at the place corresponding to the prime ideal $F_{\Fp}$ and $\varpi_{\Fp}$ be a uniformizer of the ring of integer $\CO_{F_{\Fp}}$. Let $\mathrm{H}_{\Fp^n}$ be the double coset in $\GL_2(F_{\Fp})$ with
$$\mathrm{H}_{\Fp^n}:= \GL_2(\CO_{F_{\Fp}}) \begin{pmatrix}1 &  \\  & \varpi_{\Fp^v}
\end{pmatrix} \GL_2(\CO_{F_{\Fp}}),$$
which, for $v\geqslant 1$, has measure $p^{v-1}(p+1)$ with respect to the Haar measure on $\GL_2(F_\Fp)$ assigning mass $1$ to the maximal open compact subgroup $\GL_2(\CO_{F_\Fp})$ (See \cite[Section 2.8]{sparse}). We consider the compactly supported function:
$$\mu_{\Fp^v}:= \frac{1}{p^{v/2}}\sum_{0\leqslant k\leqslant \frac{v}{2}}\mathbf{1}_{\mathrm{H}_{\Fp^{v-2k}}}.$$
Now for any $f\in \mathscr{C}^\infty(\GL_2(\BA_F))$, the Hecke operator $\Trm_{\Fp^n}$ is given by the convolution of $f$ with $\mu_{\Fp^v}$, i.e. for any $g \in \GL_2(\BA_F)$,
\begin{equation}\label{ActionHecke}
(\Trm_{\Fp^v} f)(g) = (f\star \mu_{\Fp^v})(g):= \int_{\GL_2(F_{\Fp})}f(gh)\mu_{\Fp^v}(h)dh,
\end{equation} 
and the function $h\mapsto f(gh)$ has to be understood under the natural inclusion $\GL_2(F_{\Fp})\hookrightarrow \GL_2(\BA_F)$. This definition extends to an arbitrary integral ideal $\Fh$ by multiplicativity of Hecke operators.

This abstract definition of Hecke operators simplifies a lot in the calculation when we deal with $\GL_2(\BA_F)$-invariant functionals. Indeed, consider the natural action of $\GL_2(\BA_F)$ on $\Cscr^\infty(\GL_2(\BA_F))$ by right translation and let $\ell : \Cscr^\infty(\GL_2(\BA_F))\times \Cscr^\infty(\GL_2(\BA_F))\rightarrow \BC$ be a $\GL_2(\BA_F)$-invariant bilinear functional. Then for any $f_1,f_2$ which are right $\GL_2(\CO_{F_{\Fp}})$-invariant, we have the relation
\begin{equation}\label{RelationHecke}
\ell(\Trm_{\Fp^v}f_1,f_2)= \frac{1}{p^{v/2}}\sum_{0\leqslant k\leqslant \frac{v}{2}} \gamma_{v-2k}\ell\left(\begin{pmatrix} 1 & \\ & \varpi^{v-2k}\end{pmatrix} \cdot f_1,f_2 \right),
\end{equation}
with

\begin{equation}\label{ValueGamma}
\gamma_r:= \left\{ \begin{array}{lcl}
1 & \ifm & r=0 \\ 
 & & \\
p^{r-1}(p+1) & \ifm & r\geqslant 1. 
\end{array}\right.
\end{equation}

\section{Estimations of some period integrals} \label{SectionEstimation}

Recall that $\pi_1,\pi_2$ are two unitary $\theta_i$-tempered ($i=1,2$) cuspidal automorphic representations with trivial central character and finite coprime conductor $\mathfrak {m}$ and $\mathfrak {n}$. Let $\pi_3$ be unitary $\theta_3$-tempered automorphic representation with trivial central character and finite coprime conductor $\mathfrak a$. Let $\varphi_i=\otimes_v \varphi_{i,v}\in \pi_i=\otimes_v \pi_{i,v}$ be vectors defined as follows:
We focus on the non-archimedean local fields. For $\pi_i$ ($i=1,2,3$), fix a unitary structure $\langle \cdot,\cdot\rangle_{i,v}$ on each $\pi_{i,v}$ compatible with \eqref{NormalizedInnerProduct} as in previous Section \ref{pre}. Now we need to give the choices of test vectors. Following Section 6.3 and 6.4 in \cite{hmn}, since $\Fm$ and $\Fn$ are coprime, $\pi_{1,v}$ and $\pi_{2,v}$ cannot both be ramified. If $\pi_{1,v}$ is unramified, we take $\varphi_{i,v}$ to be $L^2$-normalized newvectors for $i=2,3$ and $\varphi_{1,v}:= \pi_{1,v} \left( \begin{pmatrix} 1 & \\ & \varpi_v^s \end{pmatrix}\right) \varphi_{1,v}^0$, where $\varphi_{1,v}^0 \in \pi_{1,v}$ is the $L^2$-normalized newvector. If both $\pi_{1,v}$ and $\pi_{2,v}$ are unramified, i.e. $\crm(\pi_{1,v})=\crm(\pi_{2,v})=0$, then we simply take $\varphi_{i,v}$ to be normalized newvectors for all $i=1,2,3$. Here for large enough cardinality of the residue field, we take $s:= c(\pi_2 \otimes \pi_3)/2$, and $s:=c(\pi_2 \otimes \pi_3)/2+b$ for some non-negative absolutely bounded integer $b$ when the cardinality of the residue field is bounded. If $\pi_{2,v}$ is unramified, we take $\varphi_{i,v}$ to be normalized newvectors for $i=1,3$ and $\varphi_{2,v}:= \pi_{2,v} \left( \begin{pmatrix} 1 & \\ & \varpi_v^t \end{pmatrix}\right) \varphi_{2,v}^0$, where $\varphi_{2,v}^0 \in \pi_{2,v}$ is the normalized newvector. Here again for large enough cardinality of the residue field, we take $t:= c(\pi_1 \otimes \pi_3)/2$, and $t:=c(\pi_1 \otimes \pi_3)/2+b$ for some non-negative absolutely bounded integer $b$ when the cardinality of the residue field is bounded. It is noted that if the cardinality of the residue field is an odd prime, we may have $a \leq 1$ (See \cite{hmn}). From above choices of test vectors and for each non-archimedean local place $v$, we have the symmetry in index $i=1,2$ for the representation $\pi_{i,v}$. Hence, without loss of generality, for simplicity, we can always assume that the local representation $\pi_{1,v}$ is unramified.

In \cite{hmn}, since a uniform translation of all test vectors does not change the period integrals, they actually pick $\varphi_{1,v}$ to be normalized newvectors and $\varphi_{i,v}:= \pi_{i,v} \left( \begin{pmatrix} \varpi_v^s & \\ & 1 \end{pmatrix}\right) \varphi_{i,v}^0$, where $\varphi_{i,v}^0 \in \pi_{i,v}$ is the normalized newvector for $i=2,3$.

Let $\Fl$ be an integral ideal of $\CO_F$ which is coprime to $\Fm,\Fn,\Fa$. From the multiplicativity of the Hecke operators, without loss of generality, we simply take $\Fl$ of the form $\Fp^v$ with $\Fp\in\mathrm{Spec}(\CO_F)$ and $v\in\mathbb{N}$ and set $p$ for the norm of $\Fp$, so that $\ell=p^v$ is the norm of $\Fl$. For $0\leqslant r\leqslant v$, we write as usual 
$$\varphi_i^{\Fp^r} := \begin{pmatrix} 1 & \\ & \varpi_\Fp^r \end{pmatrix} \cdot \varphi_i.$$

\begin{rmk}\label{RemarkInfinite} 
We observe that for every finite place $v$, our local vectors $\varphi_{i,v}$ are uniquely determined, indeed there is a unique L$^2$-normalized new-vector in $\pi_v$. For the infinite place $v \vert \infty$, since we further assume that all the representations are unramified and spherical principal series, we simply pick all $\varphi_{i,v}$ to be the unique spherical vector of norm one for $i=1,2,3$, which is also the vector of minimal weight. Therefore, we make the convention that all $\ll$ involved in the following sections depend implicitly on $\pi_{i,\infty}$ for $i=1,2,3$.
\end{rmk}

\subsection{Upper and lower bounds for the local Rankin-Selberg integral and the triple product integral} \label{upper}

Before we consider the hybrid subconvexity problem for the Rankin-Selberg $L$-functions and triple product $L$-functions, we recall some results and estimations in \cite{hu}, \cite{hu1}, \cite{hu2} and \cite{hmn}.

Everything in this subsection is over non-archimedean local fields and we shall omit the subscript $v$. Let $\pi_i$ and $\pi'$ be representations of $\PGL_2$, with finite conductors.  Let $\varphi_i^0 \in \pi_i$ and $\varphi' \in \pi'$ for $i=1,2,3$ be $L^2$-normalized newvectors.

Firstly, we consider the case that $\chi$ is a character of $F^{\times}$, and the vector $\varphi_1=\varphi_{1,v} \in \pi_1= \pi(\chi, \chi^{-1},s)$ satisfies
$$ \varphi \left( \begin{pmatrix} a_1 & n \\ 0 & a_2 \end{pmatrix} g \right)= \chi(a_1)\chi^{-1}(a_2) \left \vert \frac{a_1}{a_2} \right \rvert_{F_v}^s \varphi(g).$$
When $s=\frac{1}{2}$, we simply write $\pi(\chi,\chi^{-1})=\pi(\chi,\chi^{-1},\frac{1}{2})$. In this case, we define
$$I^{RS}(\varphi_1, \varphi_2, \varphi_3):= \int_{Z(F)N \bs \GL_2(F)} W_{\varphi_2}(g) W_{\varphi_3}^{-}(g) \varphi_1(g) dg,$$
which is called the local Rankin-Selberg integral. Here $W_{\varphi}$ is the Whittaker function associated to $\varphi$ with respect to the fixed nontrivial additive character $\psi$, while $W_{\varphi}^{-}$ is for the additive character $\psi^{-}(x)=\psi(-x)$.

For general $\varphi_1 \in \pi_1$ which maynot be a principal series, we give the local triple product integral as follows:
$$I^T(\varphi_1, \varphi_2, \varphi_3):= I_v'(\varphi_{1,v} \otimes \varphi_{2,v} \otimes \varphi_{3,v}).$$

Now, we consider $\crm(\pi_2)=c>0$ and $\crm(\pi')=0$ which is unramified. Since $\crm(\pi_2)>0$ and the finite conductor of $\pi_1$ and $\pi_2$ are coprime, we automatically have $\crm(\pi_1)=0$. For the case $\crm(\pi_1)>\crm(\pi_2)=0$, we can do similarly. We will consider the following local triple product and Rankin-Selberg integral $I^T(\varphi_1, \varphi_2, \varphi_3)$, $I^T(\varphi', \varphi_2,\varphi_2)$ and $I^{RS}(\varphi', \varphi_2,\varphi_2)$. We need the upper bounds for $I^T(\varphi', \varphi_2,\varphi_2)$, $I^{RS}(\varphi', \varphi_2,\varphi_2)$ and the lower bounds for $I^T(\varphi_1, \varphi_2, \varphi_3)$ where $\varphi_{1}$ and $\varphi'$ are normalized newvectors and $\varphi_{i}:= \pi_{i} \left( \begin{pmatrix} \varpi_v^s & \\ & 1 \end{pmatrix}\right) \varphi_{i}^0$ ($\varphi_{i}^0 \in \pi_{i}$ is the normalized newvector for $i=2,3$).
They will be used later on to control the contributions from the cuspidal and the Eisenstein spectrum.

\begin{prop} \cite[Theorem 3.22]{hmn}  \label{rankintriple}
Suppose that the representation $\pi_i$ and $\pi'$ satisfy the bound $\vartheta$ towards the Ramanujan conjecture for $i=1,2,3$. 
We have the following lower bound estimation:
\begin{equation}
I_v^T(\varphi_{1,v}, \varphi_{2,v}, \varphi_{3,v}) \gg Q_v^{-\frac{1}{4}}.
\end{equation}

For the upper bounds, we have
\begin{equation}
I_v^T(\varphi_v', \varphi_{2,v}, \varphi_{2,v}) \ll_{\varepsilon} Q_v^{\varepsilon} \cdot \frac{1}{C_v(\pi_2 \otimes \pi_2)^{\frac{1}{2}}} \cdot P_v^{-\frac{1}{2}+\vartheta};
\end{equation}
and
\begin{equation}
I_v^{RS}(\varphi_v', \varphi_{2,v}, \varphi_{2,v}) \ll_{\varepsilon} Q_v^{\varepsilon} \cdot \frac{1}{C_v(\pi_2 \otimes \overline{\pi_2})^{\frac{1}{4}}} \cdot P_v^{-\frac{1}{4}+\frac{\vartheta}{2}}.
\end{equation}

\end{prop}

\begin{cor} \cite[Corollary 3.4, Remark 3.4]{BJN}  \label{square}
Suppose that $\pi_1$ is a parabolically induced representation, and $\pi_i$ satisfies the bound $\vartheta \leqslant \frac{7}{64} <\frac{1}{6}$ towards the Ramanujan conjecture for $i=1,2,3$. We further suppose that the central character of $\pi_i$ is trivial, and $\varphi_i$ is $L^2$-normalized. Then we have
$$\vert I^T(\varphi_1, \varphi_2, \varphi_3) \asymp_{\vartheta} \vert I^{RS} (\varphi_1, \varphi_2, \varphi_3) \rvert^2.$$
\end{cor}

Now we recall some results in \cite{hu}, \cite{hu1} and \cite{hu2}, which gives the double coset decomposition and integral computation.

\begin{lem}  \label{measure}
For every positive integer $c$, we have
$$\GL_2(F)= \bigsqcup_{0 \leqslant i \leqslant c} B \begin{pmatrix} 1 & 0 \\ \varpi^i & 0 \end{pmatrix} K_0(\varpi^c).$$
Here $B$ is the Borel (upper triangular) subgroup of $\GL_2$. Furthermore, if the smooth function $f$ is a $ZK_0(\varpi^c)$-invariant function, then
$$ \int_{\PGL_2(F)} f(g)dg= \sum_{0 \leqslant i \leqslant c} A_i \cdot \int_{F^{\times} \bs B(F)} f \left( b \cdot \begin{pmatrix} 1 & 0 \\ \varpi^i & 1 \end{pmatrix} \right) db.$$
Here we normalize the Haar measure on $\GL_2(F)$ such that the maximal open compact subgroup $K=\GL_2(\CO_F)$ has volume one, and $db$ is the left Haar measure on $F^{\times} \bs B(F)$ such that the subgroup $Z(\CO_F) \bs B(\CO_F)$ has volume one, and $A_0=\frac{q}{q+1}$, $A_c=\frac{1}{(q+1) \cdot q^{c-1}}$, $A_i=\frac{q-1}{(q+1)\cdot q^i}$ for $0<i<c$. Hence, $A_i \asymp q^{-i}$ are fixed constants.
\end{lem}

\subsection{Whittaker functions for newvectors over non-archimedean local fields}

Finally, we give some properties and estimations of the corresponding Whittaker functions. Here we only consider Whittaker functions over non-archimedean local fields.

\begin{lem}
Let $m \in F_v$ with $v(m)=-j<0$, and $\mu$ be a character of $\CO_F^{\times}$ with $\crm(\mu)=k>0$. Then if $j=k$, we have
$$\left \vert \int_{v(x)=0} \psi(mx) \mu^{-1}(x) d^{\times}x \right \rvert= \sqrt{q/((q-1)^2\cdot q^{k-1})}. $$
If $j \neq k$, then the corresponding integral (Gauss sum) vanishes.
\end{lem}

\begin{defn}
We define $I_{\chi,v}(x)= \chi(u)$ if $x=u \varpi^n$ for $u \in \CO_F^{\times}$. The characteristic function $I_{\chi,v}(x)$ vanishes in other cases. We say that a smooth function $f(x)$ consists of level $i$ components (with coefficients) of $L^2$-norm $h$, if we can write
$$ f(x)= \sum_{\crm(\chi)=i} \sum_{n \in \BZ} a_{\chi,v} I_{\chi,v}(x),$$
where each $\chi$ is a character of $\CO_F^{\times}$, and the norm $h=\left (\sum_{\crm(\chi)=i} \sum_{n \in \BZ} \vert a_{\chi,n} \rvert^2 \right )^{\frac{1}{2}}.$
\end{defn}

The following result is given in \cite[Proposition 2.12]{hu2}.

\begin{prop}  \label{whiproperty}
Let $\pi$ be a supercuspidal representation with $\crm(\pi)=c \geqslant 2$, or a parabolically induced representation $\pi(\mu_1, \mu_2)$ where $\crm(\mu_1)=\crm(\mu_2)=k=c/2$. Let $W$ be the $L^2$-normalized Whittaker function for a newform of $\pi$, and define
$$ W^{(i)}(a):= W \left( \begin{pmatrix} a & 0 \\ 0 & 1 \end{pmatrix} \begin{pmatrix} 1 & 0 \\ \varpi^i & 1 \end{pmatrix} \right).$$
Then we have:
\begin{enumerate}
\item $W^{(c)}(a)= I_{1,0}(a)$.
\item For $i=c-1>1$, $W^{(c-1)}(a)$ is supported only on $\CO_F^{\times}$, where it consists of level $1$ components with $L^2$-norm $\sqrt{q(q-2)/(q-1)^2}$, and also a level $0$ component with coefficient $-1/(q-1)$.
\item In general, for $0 \leqslant i<c-1$, $i \neq c/2$, $W^{(i)}(a)$ is supported only on $\{ a \in F: v(a)= \min \{ 0, 2i-c \} \}$, where it consists of level $c-i$ components with $L^2$-norm one.
\item When $i=k>1$, $W^{(c/2)}(a)$ is supported on $\CO_F$, where it consists of level $c/2$ components with $L^2$-norm one.
\item WHen $i=k=1$, $W^{(1)}(a)$ consists of a level $0$ component on $\CO_F^{\times}$ with coefficient $-1/(q-1)$, and level $1$ components on $\CO_F$ with $L^2$-norm $\sqrt{q(q-2)/(q-1)^2}$.
\end{enumerate}
\end{prop}

When the conductor of the representation $\pi$ is $0$ or $1$, we have the following very explicit formulae for the Whittaker functions.

\begin{lem}  \cite[Lemma 4.14]{hmn}  \label{unramivalue}
Suppose that $\mu_i$ are unramified (that is, $\crm(\mu_i)=0$ for $i=1,2$) and $\pi=\pi(\mu_1,\mu_2)$. Let $\varphi_0 \in \pi$ be a newform and $W_{\varphi_0}$ be its associated Whittaker function normalized so that $w_{\varphi_0}(1)=1$. Then $W_{\varphi_0}$ is invariant under the maximal open compact subgroup and
$$W_{\varphi_0}^{(0)}(\alpha)= \vert \alpha \rvert^{1/2} \times \frac{\mu_1(\varpi \alpha)-\mu_2(\varpi \alpha)}{\mu_1(\varpi)- \mu_2(\varpi)},$$
if $v(\alpha) \geqslant 0$. If $v(\alpha)<0$, we have $W_{\varphi_0}^{(0)}(\alpha)=0$.
\end{lem}

\begin{rmk}
Note that when $v(\alpha) \geqslant 0$, the numerator contains the denominator as a factor and can be cancelled. Therefore, the formula still holds when $\mu_1(\varpi)=\mu_2(\varpi)$. We also note that the above expression for $W_{\varphi_0}$ is not $L^2$-normalized, but differ only by a factor which can be controlled globally by $C(\pi)^{o(1)}$.
\end{rmk}

\begin{lem}  \cite[Lemma 4.16]{hmn}  \label{levelonevalue}
Let $\pi= \sigma(\mu \vert \cdot \rvert^{1/2}, \mu \vert \cdot \rvert^{-1/2})$ be a special unramified representation, where $\mu$ is a unramified quadratic character.

The Whittaker function associated to the newform $\varphi_0 \in \pi$ is given by
$$ W_{\varphi_0}^{(1)}(\alpha)= \mu(\alpha) \cdot \vert \alpha \rvert,$$
if $v(\alpha) \geqslant 0$. If $v(\alpha)<0$, we have $W_{\varphi_0}^{(1)}(\alpha)=0$.
$$ W_{\varphi_0}^{(0)}(\alpha)= -q^{-1} \cdot \mu(\alpha) \cdot \vert \alpha \rvert \cdot \psi(\alpha),$$
if $v(\alpha) \geqslant -1$. If $v(\alpha) \leqslant -2$, we have $W_{\varphi_0}^{(0)}(\alpha)=0$.
\end{lem}
Moreover, the value for the corresponding Whittaker function is not $L^2$-normalized, but differ only by a factor which can be controlled globally by $C(\pi)^{o(1)}$.

\section{A Symmetric Period}\label{SectionSymmetric}
Let $\pi_1,\pi_2$ and $\varphi_i\in\pi_i$ as in Section \ref{SectionEstimation}. We take $\Fq$ which is an integral ideal of $\CO_F$ and $\Fl$ which is an integral ideal of the form $\Fp^v$ with $v\in\BN$ and $\Fp \in \mathrm{Spec}(\CO_F)$ coprime with $\Fq$. Here the integral ideal $\Fq$ is the same in Section \ref{intro} and Theorem \ref{moment}. We note that integral ideal $\Fq$ is coprime to the finite coductor $C(\pi_1)$ and $C(\pi_2)$. Moreover, we recall that the integral ideal $\Fl$ is coprime to $\Fq$, $C(\pi_1)$ and $C(\pi_2)$. By the multiplicativity of the Hecke operators, without loss of generality, we write $\Fq= \Fq_1^u$ with $u \in \BN$ and $\Fq_1 \in \mathrm{Spec}(\CO_F)$ coprime with $\Fl$ and $\Fp$. We write $q, q_1, p,\ell$ for the norms of $\Fq, \Fq_1, \Fp$ and $\Fl$ respectively. We also adopt the convention that all $\ll$ involved in this section depend implicitly on the infinite datas $\varphi_{i,\infty}$ (See Remark \ref{RemarkInfinite}). Following the similar ideas in \cite{raphael} and \cite{raphael2}, we set
\begin{equation}\label{DefinitionPhi}
\Phi:=\varphi_1\varphi_2^{\Fq}.
\end{equation}
Now we consider the period
\begin{equation}\label{ThePeriod}
\CP_\Fq(\Fl,\Phi,\Phi) := \int_{X} \Trm_{\Fl}(\Phi) \overline{\Phi}=\left\langle \Trm_{\Fl}(\Phi),\Phi \right\rangle.
\end{equation}

\subsection{Expansion in the level aspect}
Since $\pi_1$ and $\pi_2$ are cuspidal, we note that $\Phi_1=\varphi_1\varphi_2^{\Fq}$ is a rapid-decay function which is invariant under the congruence subgroup $\Krm_0(\Fc[\Fm,\Fn, \Fa])$ in $\PGL_2(\BA_F)$, we apply Plancherel formula (\cite[Theorem 2.8]{raphael2}) to the well-defined inner product \eqref{ThePeriod} in the space of forms of level $\Fc [\Fm,\Fn,\Fa]$. We have the following decomposition of the considered period
\begin{equation}\label{Expansion1}
\CP_{\Fq}(\Fl,\Phi,\Phi)= \Gscr_{\Fq}(\Fl,\Phi,\Phi)+ \CC_1,
\end{equation}
where the generic part is given by
\begin{equation}\label{GenericExpansionQ}
\begin{split}
 \Gscr_{\Fq}(\Fl,\Phi,\Phi)= & \sum_{\substack{\pi \ \mathrm{cuspidal} \\ \crm(\pi)| \Fc [\Fm, \Fn, \Fa]}}\lambda_\pi(\Fl)\sum_{\psi\in\Bscr(\pi,\Fc [\Fm,\Fn, \Fa])} \left|\langle \varphi_1\varphi_2^{\Fq},\psi\rangle \right|^2 \\ 
 + & \ \sum_{\substack{\chi\in\widehat{F^\times\setminus\BA_F^{1}} \\ \crm(\chi)^2 \vert \Fc [\Fm,\Fn,\Fa]}}\int_{-\infty}^\infty \lambda_{\chi,it}(\Fl)\sum_{\psi_{it}\in\Bscr(\chi,\chi^{-1},it,\Fc [\Fm,\Fn,\Fa])}\left|\langle\varphi_1\varphi_2^{\Fq},\Erm(\psi_{it})\rangle \right|^2  \frac{dt}{4\pi}.
\end{split}
\end{equation}
Moreover, the constant term $\CC_1:= \sum_{\chi} V_F^{-1} \cdot \langle \Trm_{\Fl}(\Phi), \varphi_{\chi} \rangle \cdot \langle \varphi_{\chi}, \Phi \rangle $, where $\varphi_{\chi}(g):=\chi(\det g)$, is the one-dimensional contribution (constant term) which appears only if both $\pi_1$ and $\pi_2$ are cuspidal and there exists finite many quadratic characters $\chi$ (depending on the number field $F$) of $F^{\times} \bs \BA_F^1$ such that $\pi_1\simeq \pi_2\otimes\chi$ ($\chi=1$ if $\pi_1=\pi_2$ for example). From changing variable in the inner product $\langle \varphi_1\varphi_2^{\Fq},\varphi_{\chi} \rangle$ ($g \rightarrow gk_0$ for some $k_0$ where $k_0 \in \Krm_0(\Fb)$ is certain congruence subgroup for the finite place which stabilizes both $\varphi_1$ and $\varphi_2$), it is known that the quadratic character $\chi$ is unramified at the finite places (otherwise the constant term vanishes because of the extra term $\chi(\det(k_0))$ and a suitable choice of $k_0$). If both $\pi_1$ and $\pi_2$ are unramified at all finite places, then the problem reduced to the case in \cite{miao}. If one of $\pi_i$ ($i=1,2$) is ramified at some finite places, since the conductor of the cusp forms $\pi_1$ and $\pi_2$ are coprime, $\pi_1\simeq \pi_2\otimes\chi$ leads to the contridiction since $\chi$ is unramified at the finite places. Hence, the additional constant term $\CC_1 \equiv 0$.  Hence we conclude that
\begin{equation}\label{FirstRelation}
\CP_{\Fq}(\Fl,\Phi,\Phi)=\CG_{\Fq}(\Fl,\Phi,\Phi),
\end{equation}
where we recall that $\Phi$ is defined as \eqref{DefinitionPhi}.

\subsection{The symmetric relation} Now we follow \cite{miao} and also \cite{raphael2} to apply the symmetric relation. The symmetric relation is obtained by grouping differently the vectors $\varphi_1$ and $\varphi_2$: In the period $\CP_{\Fq}(\Fl,\Phi,\Phi)$, we first use the Hecke relation \eqref{RelationHecke} to expand the Hecke operator $\Trm_{\Fl}$. Secondly we do the same, but on the reverse way, for the translation by the matrix $\left(\begin{smallmatrix} 1 & \\ & \varpi_{\Fq_1}^u \end{smallmatrix}\right)$. Therefore, this time the Hecke operator $\Trm_{\Fq_1^m} = \Trm_{\Fq}$ appears on the dual side. By the Hecke relation \eqref{RelationHecke}, we have the following symmetric relation:

\begin{equation}  \label{Symmetry1}
q^{\frac{1}{2}} \frac{\zeta_\Fq(1)}{\zeta_\Fq(2)}\cdot  \CP_\Fq(\Fl,\Phi,\Phi) = \frac{1}{\ell^{1/2}} \cdot  \sum_{0\leqslant k\leqslant \frac{v}{2}}\gamma_{v-2k}\cdot \left(  \CP_{\Fp^{v-2k}}(\Fq,\Psi_1,\Psi_2)- \frac{1}{q_1} \cdot \CP_{\Fp^{v-2k}}(\Fq_1^{u-2},\Psi_1,\Psi_2)  \right),
\end{equation}
where we simply define $\CP_{\Fp^{v-2k}}(\Fq_1^{u-2},\Psi_1,\Psi_2):=0$ if the integral ideal $\Fq$ is squarefree, i.e. $u=1$. Moreover, 
\begin{equation}\label{Psi}
\Psi_1= \overline{\varphi}_1\varphi_1^{\Fp^{v-2k}} \ \ \mathrm{and} \ \ \Psi_2= \varphi_2\overline{\varphi}_2^{\Fp^{v-2k}}.
\end{equation}

Now we consider the period $\CP_{\Fp^{v-2k}}(\Fq,\Psi_1,\Psi_2)$ in \eqref{Symmetry1}. The period $\CP_{\Fp^{v-2k}}(\Fq_1^{u-2},\Psi_1,\Psi_2)$ on the right hand side of \eqref{Symmetry1} can be estimated in a similar way and is dominated by the period $\CP_{\Fp^{v-2k}}(\Fq,\Psi_1,\Psi_2)$. We note that the period $\CP_{\Fp^{v-2k}}(\Fq,\Psi_1,\Psi_2)$ has a similar expansion as \eqref{Expansion1}, but this time over automorphic representations of conductor dividing $\Fp^{v-2k}$. This is the phenomenon of the spectral reciprocity formula. We get a close and interesting relation between different type of $L$-functions with different spectral length. Hence, we have the following spectral decomposition:
$$\CP_{\Fp^{v-2k}}(\Fq,\Psi_1,\Psi_2)=\Gscr_{\Fp^{v-2k}}(\Fq,\Psi_1,\Psi_2)+ \mathscr{C}_2(k),$$
where $\CG_{\Fp^{v-2k}}(\Fq,\Psi_1,\Psi_2)$ is the generic part and $\mathscr{C}_2(k)$ is the constant term. 
By definition, we have

\begin{equation}\label{GenericExpansionP}
\begin{split}
  & \Gscr_{\Fp^{v-2k}}(\Fq,\Psi_1,\Psi_2) := \sum_{\substack{\pi \ \mathrm{cuspidal} \\ \crm(\pi)| \Fp^{v-2k}}}\lambda_\pi(\Fq)\sum_{\psi\in\CB(\pi,\Fp^{v-2k})} \langle \Psi_1,\psi\rangle\langle\psi, \Psi_2\rangle \\ 
 + & \; \sum_{\substack{\chi\in\widehat{F^\times\setminus\BA_F^{1}} \\ \crm(\chi)^2|\Fp^{v-2k}}}\int_{-\infty}^\infty \lambda_{\chi,it}(\Fq)\sum_{\psi_{it}\in\CB(\chi,\chi^{-1},it,\Fp^{v-2k})}\langle\Psi_1,\Erm(\psi_{it})\rangle \langle\Erm(\psi_{it}),\Psi_2\rangle \frac{dt}{4\pi}.
\end{split}
\end{equation}
Since the automorphic representation $\pi_1$ and $\pi_2$ are cuspidal, the constant term $\mathscr{C}_2(k)$ is bounded by
\begin{equation}\label{ConstantC2}
\begin{aligned}
\vert \mathscr{C}_2(k) \rvert & \leqslant \sum_{\chi} \vert \chi(\Fq) \rvert \frac{\deg(\Trm_{\Fq})\zeta_\Fq(1))}{V_F\zeta_\Fq(2)} \cdot \prod_{i=1}^2\left \vert \int_{X}\varphi_i\overline{\varphi}_i^{\Fp^{v-2k}} \varphi_{i, \chi} \right \rvert \\ 
& \leqslant \sum_{\chi}  V_{F}^{-1}q^{1/2}\frac{\zeta_\Fq(1)^2}{\zeta_\Fq(2)}\prod_{i=1}^2\left \vert \int_{X}\varphi_i\overline{\varphi}_i^{\Fp^{v-2k}} \varphi_{i, \chi} \right \rvert,
\end{aligned}
\end{equation}
where the degree of the Hecke operator $\Trm_{\Fq}$ is defined by 
\begin{equation}
\deg(\Trm_{\Fq}) := \frac{1}{q^{1/2}} \sum_{0\leqslant k\leqslant \frac{m}{2}}\gamma_{m-2k} = q^{1/2} \frac{\zeta_{F_\Fp}(1)}{\zeta_{F_\Fp}(m+1)} \leqslant q^{1/2}\zeta_{F_\Fp}(1),
\end{equation}
since $\zeta_{F_\Fp}(m+1) > 1$. Moreover, the summation is over quadratic Hecke character $\chi$ satisfying $\pi_1 \cong \pi_1 \otimes \chi$ and $\pi_2 \cong \pi_2 \otimes \chi$. We note that the cardinality of such quadratic character $\chi$ is finite (depending on the number field $F$).

For such a $\chi$, we use the identity \eqref{RelationHecke} of the Hecke relation and the definition of $\theta_2$ to have the following estimation:
\begin{equation}\label{BoundIntegral}
\left| \int_{X}\varphi_1\overline{\varphi}_2^{\Fl}\varphi_\chi\right| \leqslant \zeta_{F_{\Fp}}(1)\frac{n+1}{\ell^{1/2-\theta_2}}||\varphi_1||_{L^2}||\varphi_2||_{L^2} \ll_{\epsilon} \ell^{-1/2+\theta_2+\epsilon} \cdot ||\varphi_1||_{L^2}||\varphi_2||_{L^2}.
\end{equation}

Similarly, following \eqref{BoundIntegral}, we have
$$\left| \int_{X}\varphi_i\overline{\varphi}_i^{\Fp^{v-2k}} \varphi_{i, \chi} \right|\leqslant \zeta_{F_\Fp}(1)\frac{v-2k+1}{p^{\frac{v-2k}{2}(1-2\theta_i)}}||\varphi_i||_{L^2}^2 \Longrightarrow \Cscr_2(k)\ll_{\varepsilon, F, \pi_{1,\infty}, \pi_{2,\infty}}(mna\ell)^\varepsilon \cdot \frac{q^{1/2}}{p^{(v-2k)(1-\theta_1-\theta_2)}}.$$
Here we use the fact that $\varphi_i$ are diagonal right transpose of $L^2$-normalized vectors from our choices (See also Equation \eqref{Comparition} and \eqref{BoundAdjoint}. Hence the $L^2$-norm of $\varphi_1$ or $\varphi_2$ is bounded by $O_{\varepsilon,F,\pi_{1,\infty}, \pi_{2,\infty}}((mn)^{\varepsilon})$). Similarly, the constant term $\Cscr_3(k)$ in the period $\CP_{\Fp^{v-2k}}(\Fq_1^{u-2},\Psi_1,\Psi_2)$ can be bounded by $(mna\ell)^\varepsilon \cdot q_1^{(u-2)/2}/p^{(v-2k)(1-\theta_1-\theta_2)}$. The generic term is almost the same as \eqref{GenericExpansionP} by substituting the ideal $\Fq$ to $\Fq_1^{u-2}$.

The total constant term is obtained after summing over $0\leqslant k\leqslant n/2$ as in \eqref{Symmetry1}, i.e. 
\begin{equation}\label{GlobalConstantTerm}
\Cscr_2:= \frac{1}{\ell^{1/2}}\sum_{0\leqslant k\leqslant \frac{v}{2}}\gamma_{v-2k} \cdot \left(\Cscr_2(k)-\frac{1}{q_1} \cdot \Cscr_3(k) \right)
\end{equation}
with the following upper bound
\begin{equation}\label{BoundConstantTerm}
\mathscr{C}_2 \ll_{\varepsilon, F,\pi_{1,\infty}, \pi_{2,\infty}} (mna \ell)^\varepsilon \cdot \frac{q^{1/2}}{\ell^{1/2-\theta_1-\theta_2}}.
\end{equation}
From the above discussion, we have the following spectral reciprocity relation between the two generic parts:
\begin{equation}\label{ReciprocityRelation}
\begin{aligned}
q^{1/2}\frac{\zeta_\Fq(1)}{\zeta_\Fq(2)} \left \vert \Gscr_{\Fq}(\Fl,\Phi,\Phi) \right \rvert & \leqslant  \frac{1}{\ell ^{1/2}}\sum_{0\leqslant k\leqslant \frac{v}{2}}\gamma_{v-2k} \cdot  \vert \Gscr_{\Fp^{v-2k}}(\Fq,\Psi_1,\Psi_2) \rvert \\ &+ \Cscr_2+\frac{1}{\ell ^{1/2}}\sum_{0\leqslant k\leqslant \frac{v}{2}}\gamma_{v-2k} \cdot  \frac{1}{q_1} \left \vert \Gscr_{\Fp^{v-2k}}(\Fq_1^{u-2},\Psi_1,\Psi_2) \right \rvert.
\end{aligned}
\end{equation}

Now our main task is to bound the geometric sum $\Gscr_{\Fp^{v-2k}}(\Fq,\Psi_1,\Psi_2)$. The estimation of the geometric sum $\Gscr_{\Fp^{v-2k}}(\Fq_1^{u-2},\Psi_1,\Psi_2)$ is almost the same as $\Gscr_{\Fp^{v-2k}}(\Fq,\Psi_1,\Psi_2)$ and will give a similar bound. We need the results in Section \ref{upper} (See also \cite{hu} \cite{hu1} \cite{hu2} and \cite{hmn}) and the bound $\vert \lambda_\pi(\Fq) \rvert 
 \leqslant \tau(\Fq)q^{\vartheta}$. The main ingredients are the convexity bound for corresponding triple product $L$-functions, Proposition \ref{rankintriple} and a upper bound for certain local triple product integrals. Applying Cauchy-Schwartz inequality, it suffices to bound
 \begin{equation}\label{GenericExpansion1}
\begin{split}
S_i:= & \sum_{\substack{\pi \ \mathrm{cuspidal} \\ \crm(\pi)| \Fp^{v-2k}}}\sum_{\psi\in\CB(\pi,\Fp^{v-2k})} \left|\langle \Psi_i,\psi\rangle \right|^2 \\ 
 & + \ \sum_{\substack{\chi\in\widehat{ F^\times\setminus\BA_F^{1}} \\ \crm(\chi)^2|\Fp^{v-2k}}}\int_{-\infty}^\infty\sum_{\psi_{it}\in\CB(\chi,\chi^{-1},it,\Fp^{v-2k})}\left|\langle \Psi_i,\Erm(\psi_{it})\rangle \right|^2  \frac{dt}{4\pi},
\end{split}
\end{equation}
for $i=1,2$.

The global triple product period $S_i$ for $i=1,2$ can be bounded by the local triple product integrals from Proposition \ref{rankintriple} and the convexity bound for the triple product $L$-functions. By Proposition \ref{PropositionIntegralRepresentation}, the global triple product period is closely related to the central value of the complete triple product $L$-function. We have the following:

\begin{equation}  \label{triple8}
\begin{aligned}
& \frac{\vert \langle \Psi_i, \psi \rangle \rvert^2}{||\varphi_i||^2_{\mathrm{can}}||\varphi_i||^2_{\mathrm{can}}|| \psi ||^2_{\mathrm{can}}} \\
=& \frac{C}{8\Delta_F^{3/2}}\cdot\frac{\Lambda(\tfrac{1}{2}, \pi_i \otimes\pi_i \otimes\pi )}{ \Lambda^*(1, \pi_i,\mathrm{Ad}) \Lambda^*(1, \pi_i,\mathrm{Ad}) \Lambda^*(1, \pi ,\mathrm{Ad})}\prod_v \frac{L_v(\varphi_v)}{\langle \varphi_{i,v},\varphi_{i,v}\rangle_v \langle \varphi_{i,v},\varphi_{i,v}\rangle_v \langle \psi_v ,\psi_v \rangle_v},
\end{aligned}
\end{equation}
for $i=1,2$. Here $\psi = \otimes_v \psi_v \in \pi$ be pure tensors, $\varphi_v:= \varphi_{i,v} \otimes \varphi_{i,v} \otimes \psi_v$ and
$$L_v(\varphi_v): =L_v(\varphi_{i,v}\otimes \varphi_{i,v}\otimes \psi_v) := \zeta_{\F_v}(2)^{-2} \frac{ L(1,\pi_{i,v},\mathrm{Ad}) L(1,\pi_{i,v},\mathrm{Ad}) L(1,\pi_{v},\mathrm{Ad})}{L(\tfrac{1}{2}, \pi_{i,v}\otimes\pi_{i,v}\otimes\pi_{v})} I'_v (\varphi_{i,v}\otimes \varphi_{i,v} \otimes \psi_{v}).$$
Moreover, we have
$$I'_v(\varphi_{i,v}\otimes\varphi_{i,v}\otimes\psi_{v}) :=\int_{\PGL_2(\F_v)}\langle \pi_{i,v}(g_v)\varphi_{i,v},\varphi_{i,v}\rangle_v \langle \pi_{i,v}(g_v)\varphi_{i,v},\varphi_{i,v}\rangle_v \langle \pi_{v}(g_v)\psi_{v},\psi_{v}\rangle_v dg_v. $$
If the automorphic form $\psi$ is an Eisenstein series, we have a similar formula. We also note that it is well-known that $L(1,\pi_{i,v},\mathrm{Ad}) L(1,\pi_{i,v},\mathrm{Ad}) L(1,\pi_{v},\mathrm{Ad}) \asymp 1$ and $L(\tfrac{1}{2}, \pi_{i,v}\otimes\pi_{i,v}\otimes\pi_{v}) \asymp 1$, for example, see \cite{GJ78}.

We also need the following proposition (Proposition 11.4 in \cite{BBK}).
\begin{prop}  \label{boundtri}
Let $\pi_1$ and $\pi_2$ be unramified $\theta_i$ ($i=1,2$)-tempered principal series representations. Let $\pi_3$ be $\theta_3$-tempered of conductor $p^{\ell}$, where $p$ is the cardinality of the residue field. We further assume that $\theta_1+\theta_2+\theta_3< \frac{1}{2}$.

Let $\varphi_i \in \pi_i$, $i=1,2,3$ be the unique up-to-scalar normalized newvector. Let $W_1$ and $W_3$ be the image in the Whittaker model of $v_1$ and $v_3$.

Let $d,k \in \BZ_{\geq 0}$ be such that $\ell+d \leq k$. Then
\begin{itemize}
\item  The local triple product integral satisfies
    $$ I^T(a(\varpi^{-k}) \cdot \varphi_1, \varphi_2, a(\varpi^{-d}) \cdot \varphi_3)\ll k^4 \cdot p^{-k(1-2\theta_2)+2d(-\theta_2+\theta_1)};$$

\item   The local Rankin-Selberg integral satisfies
    $$ I^{RS}(a(\varpi^{-k}) \cdot \varphi_1, \varphi_2, a(\varpi^{-d}) \cdot \varphi_3)\ll k^2 \cdot p^{-k(1/2-\theta_2)+d(-\theta_2+\theta_1)}.$$
\end{itemize}
Here $a(\varpi^{-k})=\begin{pmatrix} \varpi^{-k} & \\ & 1 \end{pmatrix}$, which is a diagonal matrix and the action on the automorphic form is the right translation,
\end{prop}
Unconditionally, let $\theta_i=\frac{7}{64}<\frac{1}{9}$ for $i=1,2,3$. we note that we have the upper bound $I^T(a(\varpi^{-k}) \cdot \varphi_1, \varphi_2, a(\varpi^{-d}) \cdot \varphi_3)\ll k^4 \cdot p^{-5k/9}$ and $I^{RS}(a(\varpi^{-k}) \cdot \varphi_1, \varphi_2, a(\varpi^{-d}) \cdot \varphi_3)\ll k^2 \cdot p^{-5k/18}$. 

\begin{rmk} \label{special124}
Later in Section \ref{ampli}, we will see that we only need a weak form of above Proposition \ref{boundtri} in the amplification method for the proof of Theorem \ref{subconvex1}. A weaker bound $I^T(a(\varpi^{-k}) \cdot \varphi_1, \varphi_2, a(\varpi^{-d}) \cdot \varphi_3)\ll k^4 \cdot p^{-k/2}$ is enough. Moreover, for most of cases, the upper bound $I^T(a(\varpi^{-k}) \cdot \varphi_1, \varphi_2, a(\varpi^{-d}) \cdot \varphi_3)\ll 1$ will be enough for our purpose. In the step of amplification, we will see that it suffices to consider the case when $k=1,2,4$. For $k=1,2$, by explicit value of Whittaker function and Iwasawa decomposition (Lemma \ref{unramivalue} and \ref{levelonevalue}), we can estimate the local Rankin-Selberg integral explicitly. By Corollary 3.4 in \cite{BJN}, we can therefore get the upper bound for local triple product integral. We also note that for the special case $k=\ell \geq 2$ and $d=0$, the local triple product integral $ I^T(a(\varpi^{-k}) \cdot \varphi_1, \varphi_2, \varphi_3) \asymp p^{-k}$ by Theorem 1.2 and Theorem 4.1 in \cite{hu}. We consider the case $k=4$. When $\ell=3$, this is the only case which we need the weak bound $I^T(a(\varpi^{-k}) \cdot \varphi_1, \varphi_2, a(\varpi^{-d}) \cdot \varphi_3)\ll k^4 \cdot p^{-k/2}$. It suffices to prove that $I^{RS}(a(\varpi^{-k}) \cdot \varphi_1, \varphi_2, a(\varpi^{-d}) \cdot \varphi_3)\ll k^2 \cdot p^{-k/4}$.

We further assume that $k=4$, $\ell=3$, $d=1$. For the another case $k=4$, $\ell=3$, $d=0$, the proof is similar and easier. The main ingredient of the proof is the properties of certain Whittaker functions and Iwasawa decomposition (Proposition \ref{whiproperty}, Lemma \ref{unramivalue} and \ref{levelonevalue}). We sketch the proof below.

It suffices to prove that $I^{RS}(a(\varpi^{-4}) \cdot \varphi_1, \varphi_2, a(\varpi^{-1}) \cdot  \varphi_3)\ll_{\epsilon}  p^{-3/2+\epsilon}$. In order to compute the integral, we apply Lemma \ref{measure}. Since all the three representations have trivial central characters, the two diagonal matrices $a(\varpi^{-k})$ and $\begin{pmatrix} 1 & \\ & \varpi^{k} \end{pmatrix}$ can be identified as the same matrix. For fixed integers $i$ and $j$ satisfying $0 \leq i, j \leq 4$, we consider the matrix $\begin{pmatrix} 1 & \\ \varpi^i & 1 \end{pmatrix} \begin{pmatrix} 1 & \\ & \varpi^j \end{pmatrix}$.

If $i \geq j$, we have $\begin{pmatrix} 1 & \\ \varpi^i & 1 \end{pmatrix} \begin{pmatrix} 1 & \\ & \varpi^j \end{pmatrix}= \begin{pmatrix} 1 & \\ & \varpi^j \end{pmatrix} \begin{pmatrix} 1 & \\ \varpi^{i-j} & 1 \end{pmatrix}$. If $0 \leq i < j \leq 4$, we need the Iwasawa decomposition of the matrix $\begin{pmatrix} 1 & \\ \varpi^{i-j} & 1 \end{pmatrix}$. This follows from the identity
$$\begin{pmatrix} 1 & \\ \varpi^{i-j} & 1 \end{pmatrix}= \begin{pmatrix} \varpi^{j-i} & 1 \\ & \varpi^{i-j} \end{pmatrix} \begin{pmatrix} 0 & -1 \\ 1 & \varpi^{j-i} \end{pmatrix}.$$
Since $j>i$, we see that $\begin{pmatrix} 0 & -1 \\ 1 & \varpi^{j-i} \end{pmatrix} \in K$. We have the decomposition
$$ \begin{pmatrix} 0 & -1 \\ 1 & \varpi^{j-i} \end{pmatrix}= \begin{pmatrix} 1 & -1 \\ & 1 \end{pmatrix} \begin{pmatrix} 1 & 0 \\ 1 & 1 \end{pmatrix} \begin{pmatrix} 1 & -1+\varpi^{j-i} \\ 0 & 1 \end{pmatrix}.$$
Since $\begin{pmatrix} 1 & -1+\varpi^{j-i} \\ 0 & 1 \end{pmatrix} \in K_0(\varpi^4)$, we have $$\begin{pmatrix} 1 & \\ \varpi^i & 1 \end{pmatrix} \begin{pmatrix} 1 & \\ & \varpi^j \end{pmatrix} \in Z \begin{pmatrix} 1 & \varpi^{-i}-\varpi^{j-2i} \\ & 1 \end{pmatrix} \begin{pmatrix} \varpi^{j-2i} & \\ & 1 \end{pmatrix} \begin{pmatrix} 1 & 0 \\ \varpi^0 & 1 \end{pmatrix} K_0(\varpi^4).$$
Hence, for $a(t)\begin{pmatrix} 1 & \\ \varpi^i & 1 \end{pmatrix} \begin{pmatrix} 1 & \\ & \varpi^j \end{pmatrix} $, we can write
$$\begin{pmatrix} t & \\ & 1 \end{pmatrix} \begin{pmatrix} 1 & \\ \varpi^i & 1 \end{pmatrix} \begin{pmatrix} 1 & \\ & \varpi^j \end{pmatrix} \in   \begin{pmatrix} 1 & (\varpi^{-i}-\varpi^{j-2i})\cdot t \\ & 1 \end{pmatrix}  \begin{pmatrix} \varpi^{j-2i}\cdot  t & \\ & 1 \end{pmatrix} \begin{pmatrix} 1 & 0 \\ \varpi^0 & 1 \end{pmatrix} K_0(\varpi^4),$$
where we omit the center $Z$.

From above discussion, for a given Whittaker function which is right invariant by the congruence subgroup $K_0(\varpi^4)$, we have
$$ W \left( \begin{pmatrix} t & \\ & 1 \end{pmatrix} \begin{pmatrix} 1 & \\ \varpi^i & 1 \end{pmatrix} \begin{pmatrix} 1 & \\ & \varpi^j \end{pmatrix}\right)= \psi((\varpi^{-i}-\varpi^{j-2i})\cdot t) \cdot W \left( \begin{pmatrix} \varpi^{j-2i}\cdot  t & \\ & 1 \end{pmatrix} \begin{pmatrix} 1 & 0 \\ \varpi^0 & 1 \end{pmatrix} \right). $$

For above general equation for the Whittaker function, we need the special case $j=1$ and $j=4$. Now we apply Lemma \ref{measure}, Lemma \ref{unramivalue} and Proposition \ref{whiproperty} to the Whittaker function $W_2$ and $W_3$. From Lemma \ref{measure}, we consider five cases $i=0,1,2,3,4$. When $i=0$ and $i=1$, it is easy to check that the support of the Whittaker function $W_2$ and $W_3$ are disjoint. Hence, the corresponding integral equals to zero. For $i=2$, the support of the integral is $v(t)=0$. However, by Proposition \ref{whiproperty}, since $W_3$ consists of level two components with $L^2$-norm one and the level of the remaining terms are zero, hence the corresponding integral equals to zero. For $i=3$, the support of the integral is $v(t)=1$. By Proposition \ref{whiproperty}, we can have the following upper bound $1/p^3 \cdot 1/(p-1) \cdot p^{-1/2+\theta_2} \cdot p^{1/2+\theta_1} \cdot p \asymp p^{-3+\theta_1+\theta_2}< p^{-5/2}$, since we can pick $\theta_1=\theta_2=\frac{7}{64}<\frac{1}{4}$. Finally, for $i=4$, the support of the integral is $v(t)=1$. By Proposition \ref{whiproperty}, we can have the following upper bound $1/p^4 \cdot 1 \cdot p^{-1/2+\theta_2} \cdot p^{3(1/2+\theta_1)} \cdot p \asymp p^{-2+3\theta_1+\theta_2}< p^{-3/2}$, since we can pick $\theta_1=\theta_2=\frac{7}{64}<\frac{1}{8}$. Combining everything above and applying Lemma \ref{measure}, we have the upper bound for the local Rankin-Selberg integral as follows:
$$I^{RS}(a(\varpi^{-4}) \cdot \varphi_1, \varphi_2, a(\varpi^{-1}) \cdot \varphi_3) < p^{-3/2}.$$
Hence, the local triple product integral has the following upper bound
$$I^{T}(a(\varpi^{-4}) \cdot \varphi_1, \varphi_2, a(\varpi^{-1}) \cdot \varphi_3) < p^{-3}.$$
This finishes the proof and the remark.

\end{rmk}

We note that one can write an orthonormal basis in terms of linear combinations of diagonal translates of newforms, with all coefficients uniformly bounded by $k$ ($k$ is given in Proposition \ref{boundtri}), via employing the Gram–Schmidt process (See \cite[Section 7B]{nunes} and \cite[Lemma 9]{BM15} for more details and explicit coefficients). Hence, applying Proposition \ref{boundtri} and Remark \ref{special124} for each individual term, we get an nontrivial upper bound for local triple product integral when taking the summation in the orthonormal basis by spectral decomposition and Plancherel formula (for example see Equation \eqref{GenericExpansion1} and \eqref{triple8}).

When the archimedean (spectral) parameters of automorphic representation $\pi$ \eqref{GenericExpansionP} go to infinity, the corresponding triple product period integral becomes rapidly decreasing. This is because the unramified local triple product or Rankin-Selberg integrals ($I_v^T$ or $I_v^{RS}$) at archimedean places give additional Gamma factors that are rapidly decreasing by triple product formula (See Equation \eqref{triple8}, \cite{stade} and \cite{wood2}). Hence, from a weak version of the Weyl law for number fields in \cite[Section 3]{maga}, we can consider the archimedean (spectral) parameters of automorphic representation $\pi$ to be absolutely bounded.

From above discussion, especially Proposition \ref{rankintriple}, Proposition \ref{boundtri}, Remark \ref{special124}, convexity bound for triple product $L$-functions, the bound $\vert \lambda_\pi(\Fq) \rvert 
 \leqslant \tau(\Fq)q^{\theta}$, Weyl law and the Cauchy-Schwartz inequality, we obtain an upper bound for the generic terms
 $$ \Gscr_{\Fp^{v-2k}}(\Fq,\Psi_1,\Psi_2)  \ll_{\varepsilon,\F, \pi_{1,\infty}, \pi_{2,\infty}} (mna)^{\varepsilon} \cdot \left(p^{v-2k}\right)^{1+2\theta_1+2\theta_2+\varepsilon} \cdot P_{\Fq,f}^{-1/4+\frac{\theta}{2}}\cdot q^{\theta}.$$

 Finally, we can obtain the following bound:
\begin{equation}\label{FinalBound1}
\frac{1}{\ell^{1/2}} \cdot \sum_{0\leqslant k\leqslant \frac{v}{2}}\gamma_{v-2k} \cdot \Gscr_{\Fp^{v-2k}}(\Fq,\Psi_1,\Psi_2) \ll_{\varepsilon,\F, \pi_{1,\infty}, \pi_{2,\infty}} (mna\ell)^{\varepsilon} \cdot \ell^{\frac{3}{2}+2\theta_1+2\theta_2} \cdot P_{\Fq,f}^{-1/4+\frac{\theta}{2}}\cdot q^{\theta}.
\end{equation}
We also note that 
\begin{equation}
\frac{1}{\ell^{1/2}} \cdot \sum_{0\leqslant k\leqslant \frac{v}{2}}\gamma_{v-2k} \cdot \Gscr_{\Fp^{v-2k}}(\Fq_1^{u-2},\Psi_1,\Psi_2) \ll_{\varepsilon,\F, \pi_{1,\infty}, \pi_{2,\infty}} (mna\ell)^{\varepsilon} \cdot \ell^{\frac{3}{2}+2\theta_1+2\theta_2} \cdot P_{\Fq,f}^{-1/4+\frac{\theta}{2}}\cdot q^{\theta}.
\end{equation}
By the reciprocity relation (Equation \ref{ReciprocityRelation}), we deduce that
\begin{equation}  \label{Final}
\begin{aligned}
\Gscr_{\Fq}(\Fl,\Phi,\Phi) & \ll_{\varepsilon, F, \pi_{1,\infty}, \pi_{2,\infty}} (mna \ell)^{\varepsilon} \cdot \left( \ell^{-1/2+\theta_1+\theta_2}+ \ell^{\frac{3}{2}+2\theta_1+2\theta_2} \cdot P_{\Fq,f}^{-1/4+\frac{\theta}{2}}\cdot q^{-1/2+ \theta} \right) \\
& \ll_{\varepsilon, F, \pi_{1,\infty}, \pi_{2,\infty}} (mna \ell)^{\varepsilon} \cdot \left( \ell^{-1/2+\theta_1+\theta_2}+ \ell^{\frac{3}{2}+2\theta_1+2\theta_2} \cdot P_{f}^{-1/4+\frac{\theta}{2}} \right).
\end{aligned}
\end{equation}
Here we recall that the real number $\theta$ is the best exponent toward the Ramanujan-Petersson Conjecture for $\GL(2)$ over the number field $F$, we have $0 \leqslant \theta \leqslant  \frac{7}{64}.$

\subsection{Connection with the triple product}\label{Connection}
\noindent We connect in this section the expansion \eqref{GenericExpansionQ} with a first moment of the triple product $L(\tfrac{1}{2}, \pi\otimes\pi_1\otimes\pi_2)$ over automorphic representations $\pi$ of conductor dividing $\Fc \cdot [\Fm,\Fn, \Fa]$. For such a representation $\pi$, we define
\begin{equation}\label{Lscr}
\CL(\pi,\Fm,\Fn,\Fa) :=\sum_{\psi \in \CB(\pi,\Fc [\Fm,\Fn,\Fa])}\left| \langle \varphi_1\varphi_2^{\Fq},\psi\rangle \right|^2,
\end{equation}
where we recall that $\CB(\pi,\Fc [\Fm,\Fn,\Fa])$ is an orthonormal basis of the space of $\Krm_0(\Fc[\Fm,\Fn,\Fa])$-vectors in $\pi$. 
By Proposition \ref{PropositionIntegralRepresentation} and Definition \eqref{CanonicalNorm} of the canonical norm, we have
\begin{equation}\label{FactorizationLcal}
\Lscr(\pi,\Fm,\Fn,\Fa)=\frac{C}{2\Delta_F^{1/2}}f(\pi_\infty)\frac{L(\tfrac{1}{2}, \pi\otimes\pi_1\otimes\pi_2)}{\Lambda^*(1, \pi,\mathrm{Ad})} \ell(\pi,\Fm,\Fn,\Fa),
\end{equation}
where the constant $C=2 \Lambda_F(2)$. If we identify $\pi\simeq \otimes_v\pi_v$, then $\ell(\pi,\Fm,\Fn,\Fa)=\prod_{v|\Fa\Fm\Fn}\ell_v$ and the local factors $\ell_v$ are given by the summation of the local triple product integrals in Proposition \ref{PropositionIntegralRepresentation} over an orthonormal basis $\CB(\pi,\Fc [\Fm,\Fn,\Fa])$. We define the weight function $$H(\pi,\Fa,\Fm,\Fn):= \frac{\ell(\pi,\Fm,\Fn,\Fa)}{2\Delta_F^{1/2}}.$$

For certain specific weight function, if $\pi \cong \pi_3$, from the lower bound in Proposition \ref{rankintriple}, also note that $\Fm,\Fn$ are coprime and Corollary \ref{square}, we have $$ \ell(\pi_3, \Fm,\Fn,\Fa) \gg Q_f^{-\frac{1}{4}}.$$
For other automorphic representation $\pi$, by definition, it is known that $\ell_v \geqslant 0$ and the finite product $\ell(\pi, \Fm,\Fn,\Fa) \geqslant 0$.

\subsection{Archimedean function $f(\pi_\infty)$}\label{SectionInterlude}
The Archimedean function $f(\pi_\infty)$ appearing in the factorization \eqref{FactorizationLcal} is given by (See \cite[Equation (3.10)]{raphael}) 
\begin{equation}\label{Definitionfinfty}
f(\pi_\infty):= \sum_{\varphi_\infty \in \CB(\pi_\infty)}I_\infty(\varphi_\infty\otimes\varphi_{1,\infty}\otimes\varphi_{2,\infty})L(\tfrac{1}{2}, \pi_\infty\otimes\pi_{1,\infty}\otimes\pi_{2,\infty}),
\end{equation}
where the local period $I_\infty$ is defined in \eqref{DefinitionNormalizedMatrixCoefficient}. The function $f(\pi_\infty)$ is non-negative and depends on the infinite factors $\pi_{1,\infty}$ and $\pi_{2,\infty}$ and more precisely, on the choice of test vectors $\varphi_{i,\infty}\in \pi_{i,\infty}$ and the orthonormal basis $\CB(\pi_\infty)$. Since we assume that the representation $\pi_{1,\infty}$ and $\pi_{2,\infty}$ are both unramified principal series, we see that the normalized spherical vector $\varphi_{1,\infty}$ and $\varphi_{2,\infty}$ are unique. Moreover, if $f(\pi_{\infty})$ is non zero, then the representation $\pi_{\infty}$ is also unramified and $\varphi_{\infty}$ is spherical. By Corollary \ref{square}, \cite{stade} and \cite{wood2}, if $\pi_{\infty}$, $\pi_{1,\infty}$ and $\pi_{3,\infty}$ are all unramified, by picking $\varphi_{\infty}$, $\varphi_{1,\infty}$ and $\varphi_{2,\infty}$ to be normalized and spherical, we have $f(\pi_\infty) = L(\tfrac{1}{2}, \pi_\infty\otimes\pi_{1,\infty}\otimes\pi_{2,\infty})$, which is a product of Gamma functions and is rapid decay if the archimedean (spectral) parameters of automorphic representation $\pi$ go to infinity. Hence, by \cite{stade} and \cite{wood2}, we have
$$ f(\pi_{\infty}) \gg_{\pi_{1,\infty}, \pi_{2,\infty},\varepsilon} e^{-(4+\varepsilon) \cdot  \crm(\pi_{\infty})} >0.$$
Moreover, we also note that
$$ f(\pi_{\infty}) \ll_{\pi_{1,\infty}, \pi_{2,\infty},\varepsilon} e^{-(2+\varepsilon) \cdot  \crm(\pi_{\infty})}.$$

Now the estimation of Theorem \ref{moment} can be achieved from the discussion in Section \ref{SectionEstimation} and \ref{SectionSymmetric} (See also \cite[Section 4.5]{raphael2}).

\section{Proof of Theorem \ref{subconvex1}, \ref{subconvex1'}, \ref{subconvex2}}   \label{SectionSub}

We recall that $\pi_3$ is an automorphic representation of $\PGL_2(\BA_F)$ (cusp form or Eisenstein series) with finite conductor $\Fa$. 

Let $\pi_1$, $\pi_2$ be two unitary cuspidal automorphic representations satisfying the conditions in Theorem \ref{subconvex1}. We fix the test vectors $\varphi_{i} \in \pi_i$ ($i=1,2$) as in the beginning of Section \ref{SectionEstimation}.

\subsection{The amplification method} \label{ampli}
If for any positive real number $\varepsilon$, we have $P_f \ll_{\varepsilon} Q_f^{\varepsilon}$. Then there is nothing to prove and Theorem \ref{subconvex1} follows from the convexity bound for the triple product $L$-functions.

For some fixed positive real number $\delta$, we assume that $P_f \gg Q_f^{\delta}$. By Section 6.4 and Assumption 5.3 in \cite{hmn}, without loss of generality, we can further assume that $(mna)^4 \geqslant Q_f \geqslant (mna)^{1/2}$. Hence, we have $P_f \gg (mna)^{\delta/2}$.

Let $(mna)^{\delta/4048} \leqslant Q_f^{\delta/2024} \ll P_f^{1/2024}<L<P_f \leqslant Q_f^{1/2}$ be a parameter that we will choose at the end of the proof. Given $\pi$ a unitary automorphic representation of conductor dividing $\Fc   [\Fm,\Fn,\Fa]$, following \cite[Section 12]{blomerspectral} and \cite[Section 5.1]{raphael2} we choose the following amplifier
$$\CA(\pi):= \left(\sum_{\substack{\Fp \in \mathrm{Spec}(\CO_F) \\ \Nscr(\Fp)\leqslant L \\ \Fp \nmid \Fm\Fn\Fa}}\lambda_\pi(\Fp)x(\Fp)\right)^2+\left(\sum_{\substack{\Fp \in \mathrm{Spec}(\CO_F) \\ \Nscr(\Fp)\leqslant L \\ \Fp \nmid \Fm\Fn\Fa }}\lambda_\pi(\Fp^2)x(\Fp^2)\right)^2,$$
where $x(\Fl)=\mathrm{sgn}(\lambda_{\pi_3}(\Fl))$.
By Landau Prime Ideal Theorem and the Hecke relation $\lambda_{\pi_0}(\Fp)^2=\lambda_{\pi_0}(\Fp^2)+1,$ we have
\begin{equation}\label{LowerBound2}
\CA(\pi_3)\geqslant \frac{1}{2}\left(\sum_{\substack{\Fp \in \mathrm{Spec}(\CO_F) \\ \Nscr(\Fp)\leqslant L \\ \Fp \nmid \Fm \Fn \Fa }}|\lambda_{\pi_3}(\Fp)|+|\lambda_{\pi_3}(\Fp^2)|\right)^2\gg_F \frac{L^2}{(\log L)^2}.
\end{equation}
On the other hand, using the Hecke relation again, we have
\begin{equation}\label{DecompositionAmplifier}
\begin{split}
\CA(\pi)=& \; \sum_{\substack{\Fp \in \mathrm{Spec}(\CO_F) \\ \Nscr(\Fp)\leqslant L \\ \Fp \nmid \Fm \Fn \Fa }}(x(\Fp)^2+x(\Fp^2)^2)+\sum_{\substack{\Fp_1,\Fp_2  \\ \Nscr(\Fp_i)\leqslant L \\ \Fp_i \nmid \Fm \Fn \Fa }} x(\Fp_1^2)x(\Fp_2^2)\lambda_\pi(\Fp_1^2\Fp_2^2) \\ + & \; \sum_{\substack{\Fp_1,\Fp_2  \\ \Nscr(\Fp_i)\leqslant L \\ \Fp_i \nmid \Fm \Fn \Fa }} (x(\Fp_1)x(\Fp_2)+\delta_{\Fp_1=\Fp_2}x(\Fp_1^2)x(\Fp_2^2))\lambda_{\pi}(\Fp_1\Fp_2).
\end{split}
\end{equation}
Let $C$, $f(\pi_{3,\infty})$ be the quantity defined respectively in the previous Section \ref{Connection}. If $\pi_3$ is cuspidal, and by non-negativity, we have
$$C Q_f^{-1/4} \cdot \CA(\pi_3) \cdot \frac{L(\tfrac{1}{2}, \pi_1\otimes\pi_2\otimes\pi_3)}{\Lambda(1, \pi_3,\mathrm{Ad})}f(\pi_{3,\infty})\leqslant \Mscr_\CA(\pi_1,\pi_2,\Fa,\Fm,\Fn,\Fl),$$
with $\Mscr_\CA(\pi_1,\pi_2,\Fa,\Fm,\Fn,\Fl)$ as in 
\eqref{DefinitionMoment1}, but with the amplifier $\CA(\pi)$ instead of the Hecke eigenvalues in \eqref{CuspidalPart} and \eqref{ContinuousPart}. Using the lower bound \eqref{LowerBound2}, we obtain
$$\frac{L\left( \tfrac{1}{2}, \pi_1\otimes\pi_2\otimes\pi_3 \right)}{\Lambda(1, \pi_3,\mathrm{Ad})}f(\pi_{3,\infty}) \ll_{\varepsilon,\F} L^{-2+\varepsilon} \cdot Q_f^{1/4} \cdot \Mscr_\CA(\pi_1,\pi_2,\Fa,\Fm,\Fn,\Fl).$$
Now we expand the amplifier as in \eqref{DecompositionAmplifier} and apply Theorem \ref{moment} with $\Fl=1,\Fp_1\Fp_2$ or $\Fl=\Fp_1^2\Fp_2^2$ yields (Here we apply a slightly stronger version of Theorem \ref{moment} by following Remark \ref{special124} and specical integral ideal $\Fl$. Applying \cite[Theorem 4.1]{hu}, for specific integral ideal $\Fl=1,\Fp_1\Fp_2$ or $\Fl=\Fp_1^2\Fp_2^2$ which are corresponding to $k=1,2,4$ in Remark \ref{special124}, we have $ \ScM(\pi_1, \pi_2, \Fa,\Fm,\Fn, \Fl) \ll_{\pi_{1,\infty}, \pi_{2,\infty}, \F, \varepsilon} (Q_f \cdot \ell)^{\epsilon} \cdot (\ell^{\frac{3}{2}} \cdot P_{f}^{-1/4+\frac{\theta}{2}} +  \ell^{-1/2+\theta_1+\theta_2} )$, here the exponent on $\ell$ is $\frac{3}{2}$ instead of $\frac{3}{2}+2\theta_1+2\theta_2$), 
$$
\frac{L\left( \tfrac{1}{2}, \pi_1\otimes\pi_2\otimes\pi_3\right)}{\Lambda(1, \pi_3,\mathrm{Ad})}f(\pi_{3,\infty})  \ll_{\varepsilon, F, \pi_{1,\infty},\pi_{2,\infty}} Q_f^{1/4+\varepsilon} \cdot \left( L^{-1+2\theta_1+2\theta_2}+ P_{f}^{-1/4+\frac{\theta}{2}} \cdot L^6\right),
$$
Finally, choosing $L=P_f^{(1/4-\frac{\theta}{2})/(7-2\theta_1-2\theta_2)}> P_f^{1/34}$ and we get the final subconvex bound (Note that we also have $L \leqslant P_f^{1/28}$):
\begin{equation}\label{HybridBound}
\frac{L\left( \tfrac{1}{2}, \pi_1\otimes\pi_2\otimes\pi_3 \right)}{\Lambda(1, \pi_3,\mathrm{Ad})}f(\pi_{3,\infty})  \ll_{\varepsilon, F, \pi_{1,\infty},\pi_{2,\infty}} Q_f^{1/4+\varepsilon} \cdot P_f^{-(\frac{1}{4}-\frac{\theta}{2})(1-2\theta_1-2\theta_2)/(7-2\theta_1-2\theta_2)}.
\end{equation}
Using \eqref{BoundAdjoint} for the adjoint $L$-function at $s=1$ and Section \ref{SectionInterlude}, equation \eqref{HybridBound} transforms into
$$L\left( \tfrac{1}{2}, \pi_1 \otimes\pi_2\otimes\pi_3 \right) \ll_{\varepsilon, F,\pi_{1,\infty},\pi_{2,\infty},\pi_{3,\infty}} Q_f^{1/4+\varepsilon} \cdot P_f^{-(\frac{1}{4}-\frac{\theta}{2})(1-2\theta_1-2\theta_2)/(7-2\theta_1-2\theta_2)},$$
which gives the desired subconvexity bound in Theorem \ref{subconvex1} for the cuspidal part.

If $\pi_3$ is not cuspidal, i.e. an Eisenstein series, the proof of Theorem \ref{subconvex1} for the Eisenstein part is almost the same as above. Instead of the cuspidal distribution and its postivity (See Equation \ref{CuspidalPart}), we need the continuous distribution and its non-negativity (See Equation \ref{ContinuousPart}). From the discussion between Equation \ref{ContinuousPart} and Equation \ref{DefinitionMoment1}, we see that the only obstacle in deducing a subconvex bound is as follows: When the variable $t$ attaches to $0$ and $\chi$ is a quadratic character, the quotient $L(\tfrac{1}{2}+it, \pi_1\otimes\pi_2\otimes\omega)L(\tfrac{1}{2}-it, \pi_1\otimes\pi_2\otimes\omegabar)/{\Lambda^*(1, \pi_\omega(it),\mathrm{Ad})}= \vert L(\tfrac{1}{2}+it, \pi_1\otimes\pi_2\otimes\omega) \rvert^2/ {\Lambda^*(1, \pi_\omega(it),\mathrm{Ad})}$ has a zero of order two at $t=0$. One can overcome this barrier by an application of Holder's inequality, as in \cite[Section 4]{blo}.

For Theorem \ref{subconvex1'} and Theorem \ref{subconvex2}, they are actually corollaries of Theorem \ref{subconvex1}, since we note that when $C(\pi)=p^c$ and $C(\pi \otimes \overline{\pi})=C(\pi \otimes \pi)=p^{d}$, then $d$ is an even number and satisfies $d \leqslant c+1$ (For more details and when $d=c+1$, see \cite[Proposition 2.5]{nps}).


\begin{thebibliography}{xinchenmiao}


\bibitem[AK18]{andersen}
Nickolas Andersen and Eren~Mehmet Kiral.
\newblock Level reciprocity in the twisted second moment of {R}ankin-{S}elberg
  {$L$}-functions.
\newblock {\em Mathematika}, 64(3):770--784, 2018.


\bibitem[Blo12]{blo}
Valentin Blomer.
\newblock Subconvexity for twisted $L$-functions on $\GL(3)$.
\newblock {\em American Journal of Math.}, 134(5):1385-1421, 2012.


\bibitem[BB11]{ramanujan}
Valentin Blomer and Farrell Brumley.
\newblock On the {R}amanujan conjecture over number fields.
\newblock {\em Ann. of Math. (2)}, 174(1):581--605, 2011.

\bibitem[BBK22]{BBK}
Valentin Blomer, Farrell Brumley and Ilya Khayutin.
\newblock The mixing conjecture under GRH.
\newblock {\em arXiv preprint arXiv: 2212.06280}, 2022.

\bibitem[BHKM20]{moto}
Valentin Blomer, Peter Humphries, Rizwanur Khan, and Micah Milinovich.
\newblock Motohashi's fourth moment identity for non-archimedean test functions
  and applications.
\newblock {\em arXiv preprint arXiv: 1902.07042}, 2019.

\bibitem[BJN24]{BJN}
Valentin Blomer, Subhajit Jana and Paul D. Nelson.
\newblock Local integral transforms and global spectral decomposition.
\newblock {\em arXiv preprint arXiv: 2404.10692}, 2024.

\bibitem[BK19]{blomerspectral}
Valentin Blomer and Rizwanur Khan.
\newblock Twisted moments of {$L$}-functions and spectral reciprocity.
\newblock {\em Duke Math. J.}, 168(6):1109--1177, 2019.

\bibitem[BLM19]{blomer4}
Valentin Blomer, Xiaoqing Li, and Stephen~D. Miller.
\newblock A spectral reciprocity formula and non-vanishing for {$L$}-functions
  on {${\rm GL}(4)\times {\rm GL}(2)$}.
\newblock {\em J. Number Theory}, 205:1--43, 2019.

\bibitem[BM15]{BM15}
Valentin Blomer and Djordje Milićević.
\newblock The Second Moment of Twisted Modular $L$-Functions.
\newblock {\em  Geom. Funct. Anal.}, 25(2) 453-516, 2015.



\bibitem[Bum97]{bump}
Daniel Bump.
\newblock {\em Automorphic forms and representations}, volume~55 of {\em
  Cambridge Studies in Advanced Mathematics}.
\newblock Cambridge University Press, Cambridge, 1997.


\bibitem[Gel75]{adele}
Stephen~S. Gelbart.
\newblock {\em Automorphic forms on ad\`ele groups}.
\newblock Princeton University Press, Princeton, N.J.; University of Tokyo
  Press, Tokyo, 1975.
\newblock Annals of Mathematics Studies, No. 83.

\bibitem[GJ78]{GJ78}
Stephen Gelbart and Herv\'e Jacquet.
\newblock A relation between automorphic representations of $\GL(2)$ and $\GL(3)$.
\newblock Annales Scientifiques de I'ENS Volume 11 Issue 4, page 471-542.


\bibitem[GJ79]{analytic}
Stephen Gelbart and Herv\'e Jacquet.
\newblock Forms of {${\rm GL}(2)$}\ from the analytic point of view.
\newblock In {\em Automorphic forms, representations and {$L$}-functions
  ({P}roc. {S}ympos. {P}ure {M}ath., {O}regon {S}tate {U}niv., {C}orvallis,
  {O}re., 1977), {P}art 1}, Proc. Sympos. Pure Math., XXXIII, pages 213--251.
  Amer. Math. Soc., Providence, R.I., 1979.







\bibitem[HM06]{michel2}
Gergely Harcos and Philippe Michel.
\newblock The subconvexity problem for {R}ankin-{S}elberg {$L$}-functions and
  equidistribution of {H}eegner points. {II}.
\newblock {\em Invent. Math.}, 163(3):581--655, 2006.

\bibitem[HL94]{adjoint}
Jeffrey Hoffstein and Paul Lockhart.
\newblock Coefficients of {M}aass forms and the {S}iegel zero.
\newblock {\em Ann. of Math. (2)}, 140(1):161--181, 1994.
\newblock With an appendix by Dorian Goldfeld, Hoffstein and Daniel Lieman.

\bibitem[HM13]{holo}
Roman Holowinsky and Ritabrata Munshi.
\newblock Level aspect subconvexity for {R}ankin-{S}elberg {$L$}-functions.
\newblock In {\em Automorphic representations and {$L$}-functions}, volume~22
  of {\em Tata Inst. Fundam. Res. Stud. Math.}, pages 311--334. Tata Inst.
  Fund. Res., Mumbai, 2013.

\bibitem[Hu17]{hu}
Yueke Hu.
\newblock Triple product formula and the subconvexity bound of triple product
  {$L$}-function in level aspect.
\newblock {\em Amer. J. Math.}, 139(1):215--259, 2017.

\bibitem[Hu18]{hu1}
Yueke Hu.
\newblock Triple product formula and mass equidistribution on modular curves of level $N$.
\newblock {\em Int. Math. Res. Not.}, 9, 2899--2943, 2018.

\bibitem[Hu20]{hu2}
Yueke Hu.
\newblock Mass equidistribution on the torus in the depth aspect.
\newblock {\em Algebra and Number Theory}, 14 (4), 2020.

\bibitem[HMN23]{hmn}
Yueke Hu, Philippe Michel, and Paul Nelson.
\newblock The subconvexity problem for Rankin-Selberg and triple product $L$-functions.
\newblock {\em arXiv: 2207.14449v2}.


\bibitem[Ich08]{ichino}
Atsushi Ichino.
\newblock Trilinear forms and the central values of triple product
  {$L$}-functions.
\newblock {\em Duke Math. J.}, 145(2):281--307, 2008.

\bibitem[JS81]{classification}
Herve~Jacquet and Joseph~A. Shalika.
\newblock On {E}uler products and the classification of automorphic
  representations. {I}.
\newblock {\em Amer. J. Math.}, 103(3):499--558, 1981.

\bibitem[Kha12]{simult}
Rizwanur Khan.
\newblock Simultaneous non-vanishing of {$\GL(3)\times \GL(2)$} and {$\GL(2)$}
  {$L$}-functions.
\newblock {\em Math. Proc. Cambridge Philos. Soc.}, 152(3):535--553, 2012.


\bibitem[KMV02]{michel0}
Emmanuel~Kowalski, Philippe~Michel, and Jeffrey~VanderKam.
\newblock Rankin-{S}elberg {$L$}-functions in the level aspect.
\newblock {\em Duke Math. J.}, 114(1):123--191, 2002.

\bibitem[Maga17]{maga}
Peter Maga.
\newblock Subconvexity for twisted $L$-functions over number fields via shifted convolution sums.
\newblock {\em Acta Math. Hungar.}, 151 (1) 232-257, 2017.


\bibitem[Miao24]{miao}
Xinchen Miao.
\newblock Spectral reciprocity for the first moment of triple product $L$-functions and applications.
\newblock preprint, 2024.

\bibitem[Mic04]{michel1}
Philippe~Michel.
\newblock The subconvexity problem for {R}ankin-{S}elberg {$L$}-functions and
  equidistribution of {H}eegner points.
\newblock {\em Ann. of Math. (2)}, 160(1):185--236, 2004.

\bibitem[Mic22]{mic}
Philippe Michel.
\newblock Recent progresses on the subconvexity problem.
\newblock {\em Seminar Bourbaki}, no.1190, (74), 2021-2022.


\bibitem[MV10]{subconvexity}
Philippe Michel and Akshay Venkatesh.
\newblock The subconvexity problem for {${\rm GL}_2$}.
\newblock {\em Publ. Math. Inst. Hautes \'{E}tudes Sci.}, (111):171--271, 2010.

\bibitem[Nel19]{nelson}
Paul~D. Nelson.
\newblock Subconvex equidistribution of cusp forms: reduction to {E}isenstein observables.
\newblock {\em Duke Math. J.}, 168(9):1665--1722, 2019.


\bibitem[NPS14]{nps}
Paul D. Nelson, Ameya Pitale and Abhishek Saha.
\newblock Bounds for Rankin–Selberg integrals and quantum unique ergodicity for powerful levels.
\newblock {\em J. Amer. Math. Soc.}, 27 (2014), 147-191.




\bibitem[Nun23]{nunes}
Ramon Nunes.
\newblock Spectral reciprocity via integral representations.
\newblock {\em 	Alg. Number Th.}, 17(8):1381-1409, 2023.

\bibitem[Sta02]{stade}
Eric Stade.
\newblock Archimedean $L$-factors on $\GL(n) \times \GL(n)$ and Generalized Barnes Integrals.
\newblock {\em Israel Journal of Math.}, 127(1), 201-219, 2002.


\bibitem[Ven10]{sparse}
Akshay Venkatesh.
\newblock Sparse equidistribution problems, period bounds and subconvexity.
\newblock {\em Ann. of Math. (2)}, 172(2):989--1094, 2010.

\bibitem[Wood11]{wood}
Michael~C. Woodbury.
\newblock  Explicit trilinear forms and subconvexity of the triple product $L$-function.
\newblock ProQuest LLC, Ann Arbor, MI, 2011.
\newblock Thesis (Ph.D.)--The University of Wisconsin - Madison.

\bibitem[Wood12]{wood3}
Michael~C. Woodbury.
\newblock Trilinear forms and subconvexity of the triple product $L$-functions.
\newblock submitted, 2012.


\bibitem[Wood17]{wood2}
Michael~C. Woodbury.
\newblock  On the triple product formula: real local calculations.
\newblock preprint 2017.


\bibitem[Wu14]{han2}
Han Wu.
\newblock Burgess-like subconvex bounds for {$\text{GL}_2\times\text{GL}_1$}.
\newblock {\em Geom. Funct. Anal.}, 24(3):968--1036, 2014.





\bibitem[Zac19]{raphael}
Raphaël Zacharias.
\newblock Periods and reciprocity {I}.
\newblock {\em International Mathematics Research Notices},
\newblock rnz100.

\bibitem[Zac20]{raphael2}
Raphaël Zacharias.
\newblock Periods and reciprocity {II}.
\newblock {\em arXiv: 1912.01512v2}.



\end{thebibliography}
\end{document}